\documentclass[3p, number,sort&compress,times,fleqn]{elsarticle}

\usepackage{amsmath}
\usepackage{amssymb}
\usepackage{mathtools}
\usepackage[colorlinks = true,
linkcolor = blue,
urlcolor  = blue,
citecolor = blue,
anchorcolor = blue]{hyperref}

\usepackage{subfig}
\usepackage{graphicx}
\usepackage[labelfont=bf, justification=justified]{caption}
\usepackage{booktabs}
\usepackage{threeparttable}
\usepackage{float}

\usepackage{./style/csml}

\graphicspath{{./figs/}}
\usepackage{url}

\begin{document}
	
\begin{frontmatter}
	
\title{An optimally convergent smooth blended B-spline construction for semi-structured quadrilateral and hexahedral meshes}

\author[1]{Kim Jie Koh}
\author[2]{Deepesh Toshniwal}
\author[1]{Fehmi Cirak\corref{cor1}}
\ead{f.cirak@eng.cam.ac.uk}

\cortext[cor1]{Corresponding author}

\address[1]{Department of Engineering, University of Cambridge, Cambridge, CB2 1PZ, UK }
\address[2]{Delft Institute of Applied Mathematics, TU Delft, 628 XE Delft, The Netherlands}

\begin{abstract}
Easy to construct and optimally convergent generalisations of B-splines to unstructured meshes are essential for the application of isogeometric analysis to domains with non-trivial topologies. Nonetheless, especially for hexahedral meshes, the construction of smooth and optimally convergent isogeometric analysis basis functions is still an open question. We introduce a simple partition of unity construction that yields smooth blended B-splines, referred to as SB-splines, on semi-structured quadrilateral and hexahedral meshes, namely on mostly structured meshes with a few sufficiently separated unstructured regions. To this end, we first define the mixed smoothness B-splines that are $C^0$ continuous in the unstructured regions of the mesh but have higher smoothness everywhere else. Subsequently, the SB-splines are obtained by smoothly blending in the physical space the mixed smoothness B-splines with Bernstein bases of equal degree. One of the key novelties of our approach is that the required smooth weight functions are assembled from the available smooth B-splines on the unstructured mesh. The SB-splines are globally smooth, non-negative, have no breakpoints within the elements and reduce to conventional B-splines away from the unstructured regions of the mesh. Although we consider only quadratic mixed smoothness B-splines in this paper, the construction generalises to arbitrary degrees. We demonstrate the excellent performance of SB-splines studying Poisson and biharmonic problems on semi-structured quadrilateral and hexahedral meshes, and numerically establishing their optimal convergence in one and two dimensions.
\end{abstract}
	

\end{frontmatter}


%
\section{Introduction \label{sec:introduction}}
\begin{figure}[!b]
	\centering
	\subfloat[][Hexahedral mesh \label{fig:introSphereControl}] {
		\includegraphics[width=0.3 \textwidth]{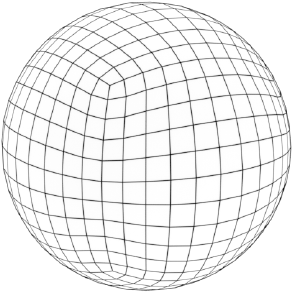} }
	\hfill
	\subfloat[][B\'ezier mesh \label{fig:introSphereParam}] {
		\includegraphics[width=0.3 \textwidth]{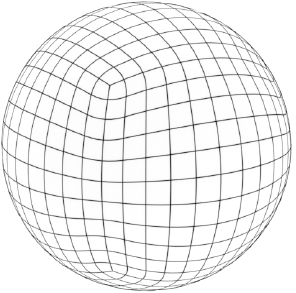} }
	\hfill
	\subfloat[][Extraordinary edges and vertices  \label{fig:introSphereExtraordinary}] {
		\includegraphics[width=0.3 \textwidth]{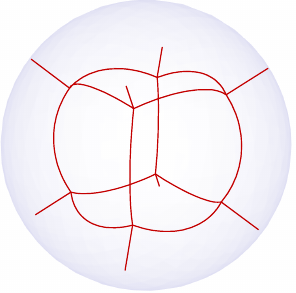} }
	\caption{ Spherical domain discretised with a hexahedral mesh containing extraordinary edges and vertices. The B\'ezier mesh represents the physical domain.\label{fig:introSphere}}
\end{figure}

%
The smoothness of spline basis functions is vital in the isogeometric analysis of problems with higher-order partial differential equations. For instance, gradient-theories of elasticity and plasticity~\cite{fischer2011isogeometric,rudraraju2014three,niiranen2016variational,de2016gradient,codony2019immersed}, phase-field modelling of sharp interfaces~\cite{gomez2008isogeometric,dede2012isogeometric,liu2013isogeometric} and Kirchhoff-Love type plate and shell models~\cite{Cirak:2000aa,kiendl2009isogeometric,benson2011large,bartezzaghi2015isogeometric} and their extensions~\cite{Long:2012aa, echter2013hierarchic} all lead to  higher-order partial differential equations. Since its inception, isogeometric analysis brought about a revival of such theories mainly because of the ease of discretising higher-order partial differential equations using smooth spline basis functions.  In particular, smooth basis functions avoid the introduction of (non-physical) extra degrees of freedom and promise a better integration with common computer-aided design representations. However, standard spline basis functions, including B-splines, NURBS and box-splines, are defined only on structured meshes and must be suitably extended for domains with non-trivial topology. For instance, multivariate B-splines are defined only on structured quadrilateral and hexahedral meshes in 2D and 3D, respectively. Most industrial complex geometries cannot be parametrised with a structured mesh so that a limited number of singularities on the surface or inside the volume must be introduced~\cite{murdoch1997spatial,tarini2004polycube,nieser2011cubecover,shepherd2008hexahedral,li2012all,zhang2018geometric,bracci2019hexalab,zhang2020octahedral,livesu2020loopycuts}. These singularities manifest themselves as extraordinary vertices and edges in the mesh, see Figure~\ref{fig:introSphere}. For a hexahedral mesh, an interior vertex is extraordinary if it is not incident to $8$ hexahedra, and an interior edge is extraordinary  if it is not incident to $4$ hexahedra. Similarly, an interior vertex is extraordinary for a quadrilateral mesh if it is not incident to $4$. The construction of smooth splines which generalise or extend B-splines to unstructured meshes is currently a very active area of research in isogeometric analysis.

In computer-aided design numerous constructions have been proposed to deal with extraordinary vertices in a surface mesh, including geometrically~$G^k$ and parametrically~$C^k$ continuous constructions~\cite{derose1990necessary,reif1998turbs,scott2014isogeometric,nguyen2016c1,collin2016analysis,toshniwal2017smooth,toshniwal2017multi,kapl2017isogeometric,kapl2018construction,karvciauskas2021multi}, subdivision surfaces~\cite{Doo:1978aa,Catmull:1978aa,Peters:2008aa,wei2015truncated,wei2021tuned,zhang2018subdivision,ma2019subdivision}, macro-elements~\cite{clough1965finite,powell1977piecewise,lai2007spline} and manifold constructions~\cite{grimm1995modeling,della2008construction,ying2004simple,tosun2011manifold,majeedCirak:2016,zhang2020manifold,zhang2019manifold}. There is,  however, a very limited number of constructions for volume meshes, including \cite{bajaj2002subdivision,chang2002new,xie2020interpolatory,reif2019old,peters2020refinable,wei2018blended}; most likely because conventional computer-aided design representations do not require a volume parametrisation. As widely reported, most constructions from computer-aided design do not lead in isogeometric analysis to optimally convergent finite elements, especially when applied to higher-order partial differential equations, see the discussion in~\cite{toshniwal2017smooth}. There are, however, constructions for unstructured quadrilateral meshes, including~\cite{toshniwal2017multi,kapl2018construction,majeedCirak:2016, zhang2020manifold,wei2021tuned}, which yield optimal or nearly optimal convergence rates. In contrast, there are no B-spline based  optimally convergent smooth constructions for unstructured hexahedral meshes. Currently,  optimality is achieved by either reducing continuity to~$C^0$ around extraordinary features~\cite{wei2018blended,schneider2019poly,schneider2021isogeometric}, combining B-splines with meshless approximants~\cite{wang2014consistently} or resorting to non-standard spline definitions~\cite{febrianto2021mollified}. The first approach is not suitable for higher-order partial differential equations whereas the latter approaches lead to schemes that are usually computationally very costly.
\begin{table} 
	\centering
	\begin{threeparttable}
		\caption{Summary of the terminology and symbols used in this paper. Note that this terminology, in particular the usage of \emph{blending} or \emph{blended} is not uniform across the isogeometric analysis and geometric design literature.  \label{tab:terminology}}
	\begin{tabular}{ll}
		\toprule
		Terminology & Definition
		\\ \midrule
		Mixed B-splines, $\vec B(\vec x)$ & B-splines of mixed smoothness.  \\
		Blending & Combination of different basis functions using weight  functions.  \\
		SB-splines, $\vec N(\vec x)$ & Smooth blended splines obtained by the proposed construction. \\ \bottomrule
	\end{tabular}
	\end{threeparttable}
\end{table}

In this paper, we derive a computationally efficient, easy to construct and optimally convergent extension of B-splines to unstructured quadrilateral and hexahedral meshes. We dub the new basis functions \emph{SB-splines}; see Table~\ref{tab:terminology} for the terminology used throughout this paper. Although we consider only quadratic B-splines, the presented ideas should carry over to arbitrary degrees. To begin with, we determine on the given unstructured mesh a set of B-splines of mixed smoothness following the construction for quadrilateral meshes presented in Toshniwal~\cite{toshniwal2022quadratic}. The mixed B-splines are~$C^0$ continuous around extraordinary features, i.e. extraordinary vertices (in 2D and 3D) and extraordinary edges (in 3D), but are~$C^1$ smooth everywhere else. Subsequently, we use the partition of unity method of Melenk and  Babu{\v{s}}ka~\cite{melenk1996partition} to blend the mixed B-splines~$\vec B(\vec x)$ with tensor-product Bernstein basis functions~\mbox{$\vec Q(\vec x)$} of equal degree.  To this end, a set of smooth weight, or partition of unity, functions~$w^B(\vec x)$ and~$w^Q(\vec x)$ are defined to blend both types of basis functions. A key novelty of our approach is that the blending function~\mbox{$w^B(\vec x)$} is assembled from the mixed B-splines on the unstructured mesh, by excluding the~$C^0$ ones, so that~\mbox{$w^Q(\vec x) \coloneqq  1- w^B(\vec x)$}. Consequently, the weight functions~$w^B(\vec x)$ and~$w^Q(\vec x)$ are~$C^1$ smooth and have their breakpoints at the element boundaries.  On quadrilateral meshes with extraordinary vertices the SB-splines are simply given by the weighted basis functions ~\mbox{$w^B(\vec x) \vec B(\vec x)$} and~\mbox{$w^Q(\vec x) \vec Q(\vec x)$}.  In hexahedral meshes, the extraordinary edges and vertices usually form a connected network as illustrated in Figure~\ref{fig:introSphereExtraordinary}; see also relevant work on hexahedral meshing~\cite{murdoch1997spatial,tarini2004polycube,nieser2011cubecover,shepherd2008hexahedral,li2012all,zhang2018geometric,bracci2019hexalab,zhang2020octahedral,livesu2020loopycuts}. That is, the weight function \mbox{$w^Q(\vec x) =  1- w^B(\vec x)$} has a support over the entire network and is decomposed as \mbox{$ w^Q(\vec x) = \sum_i \sum_j w_{i,j}^P(\vec x) + \sum_j w_j^J(\vec x)$}  into locally supported weight functions. In brief, the support of one weight function~\mbox{$ w_j^J(\vec x)$} covers a region where more than two extraordinary edges meet, and the supports of the weight functions~\mbox{$ w_{i,j}^P(\vec x) $} are restricted to the remaining regions along the connecting extraordinary edge chains. We note for blending instead of the Bernstein basis~$\vec Q( \vec x)$ a different tensor-product basis or a triangular Bernstein-B\'ezier basis can be considered. Similarly, in principle, differently constructed mixed smoothness B-splines~$\vec B(\vec x)$ could be used, including the enhanced smoothness B-splines presented in Buchegger and J\"uttler  \cite{buchegger2016adaptively}, which are conceptually similar to the chosen mixed B-splines. However, the enhanced smoothness B-splines use nested refinement (leading to a more involved implementation because of the mixed smoothness) and exhibit optimal convergence only for extraordinary vertices with a valence $3$ using uniform refinement.

The proposed approach uses, like the manifold-based constructions~\cite{ ying2004simple,tosun2011manifold,majeedCirak:2016,zhang2019manifold, zhang2020manifold} and their variations~\cite{zorin2006constructing,levin2006modified,pla2006n,antonelli2013subdivision}, the partition of unity method to smoothly blend mixed B-splines with $C^\infty$ continuous Bernstein basis functions. Unlike manifold-based constructions, the two types of basis functions are blended in the Euclidean ambient space, i.e.~$\mathbb R^2$ in 2D and~$\mathbb R^3$ in 3D. Thus, for blending we do not use an atlas consisting of charts and smooth transitions maps. Although it is easy to devise smooth transition maps for 2-manifolds, e.g. using conformal or characteristic maps~\cite{ying2004simple,Peters:2008aa}, it is not clear how to construct them in~$\mathbb R^3$. Circumventing the need for such a smooth atlas yields a conceptually and implementation-wise simpler approach and smooth basis functions with appealing properties. The SB-splines are obtained by blending polynomials defined either on the parameter or ambient space and, hence, can be integrated very efficiently using standard Gauss-Legendre quadrature. Furthermore, the weight functions~$w^B(\vec x)$ and~$w^Q(\vec x)$ have minimal polynomial degree considering that they are assembled from B-splines and their complements to one.

The outline of this paper is as follows. To begin with, we briefly discuss in Section~\ref{sec:oneD} the construction of SB-splines in 1D to introduce the key ideas and terminology used throughout the paper. Subsequently,  we consider in Section~\ref{sec:two} the construction of SB-splines on unstructured quadrilateral meshes with extraordinary vertices. We first review the mixed B-splines  in Section~\ref{sec:C0blended} and then discuss the construction of SB-splines, in particular the weight functions, in Section~\ref{sec:C1blended}. As explained in Section~\ref{sec:C1blended}, it is straightforward to derive closed form expressions of the new basis functions for use in existing isogeometric analysis implementations. In Section~\ref{sec:three},  we consider the construction of SB-splines on unstructured hexahedral meshes with extraordinary edges and vertices. After discussing the extension of mixed B-splines  to hexahedral meshes in Section~\ref{sec:C0blendedHex},  we first introduce the notion of an \emph{extraordinary prism} and ~\emph{extraordinary joint} and discuss how to construct the respective weight and smooth basis functions in Section~\ref{sec:C1blendedHex}. Finally, we introduce in Section~\ref{sec:examples} several Poisson and biharmonic examples to confirm the convergence of the SB-splines. We study in particular the influence of the number of quadrature points and the valence of the extraordinary vertices on finite element convergence in 2D. Although we have not investigated the finite element convergence in 3D, we demonstrate the global $C^1$ continuity and excellent performance of the finite element solution on the same spherical domain as shown in Figure \ref{fig:introSphere}. The paper is supplemented by four appendices that provide a proof of linear independence and discuss aspects of finite element discretisation, mesh refinement and an illustration of the treatment of arbitrary hexahedral meshes.

\section{One-dimensional SB-splines \label{sec:oneD}}
%
The proposed construction is best illustrated in the one-dimensional setting.  Given is a domain~\mbox{$\Omega \subset \mathbb{R}$} with the parametrisation
\begin{equation} \label{eq:oneDMapping}
x(\xi) = \sum_{i = 1}^{n_{B}} B_i(\xi) x_i  \quad~\text{with } \xi \in \hat\Omega:=[0,1] \,,
\end{equation}
where~$B_i(\xi)$ are the univariate B-splines of degree~$p_B \ge 2$ and~$x_i$ are the coordinates of the control points. As mentioned earlier, in higher dimensions we will focus solely on the case~$p_B=2$.  The B-splines~$B_i(\xi)$ are defined on the parametric domain~$\hat \Omega$ with the parametric coordinate~$\xi$. For the sake of illustration, they are chosen to be $C^0$ continuous at the break point \mbox{$\xi_{\text{ep}}$} and $C^{p_B-1}$  continuous at every other break point.  Hence, the point with the coordinate~$x_\text{ep} \coloneqq  x(\xi_\text{ep}) $ is an extraordinary point.

 Assuming that the parametrisation~$x(\xi)$ is, as usual, bijective,  the push-forward of the basis functions on the physical domain of interest~$\Omega$ are given by
\begin{equation} \label{eq:oneDMappingBSplines}
	B_i ( x ) =  B_i(\xi)  \circ x(\xi)^{-1}
\end{equation}
For later reference, the basis functions on the physical domain~$\Omega$ are collected in the array
\begin{equation} \label{eq:oneDUnivariateBSplines}
\vec B(x) = 
\begin{pmatrix}
	B_1(x) & \dotsc & B_{n_B}(x) 
\end{pmatrix}^\trans  \, .
\end{equation}
In the neighbourhood of the extraordinary point~$x_\text{ep}$, we aim to blend the B-splines~$\vec B(x)$ with a second polynomial basis defined only over the the blending domain~$ \Omega^Q \subset \Omega$. The second basis is, without loss of generality, throughout this paper a Bernstein basis with the basis functions
\begin{equation} \label{eq:oneDEnrichmentFunctions}
\vec Q(x) = \begin{pmatrix}
	Q_1(x) & \dotsc & Q_{n_Q}(x) 	\end{pmatrix}^\trans   \, , 
\end{equation}
where~$n_Q = p_Q +1$ and~$p_Q$ is the polynomial degree of the Bernstein basis. 

For blending together the two sets of basis functions, we choose a weight function~$w^Q(x)$ with~\mbox{$\supp w^Q(x) = \Omega^Q$}  and its complement to one~$w^B(x) $, i.e.,
\begin{equation}
	w^Q(x) + w^B(x) \equiv 1 \quad \forall x \in \Omega \,.
\end{equation}
Such weight functions can be chosen in many different ways. In the proposed construction, we assemble the weight function~$w^B(x)$ from the B-splines~$B_i (x)$. In particular, with the index set~\mbox{$\set I = \{ i \mid B_i (x) \text{ is at most } C^0 \text{ at } x_\text{ep}  \}$} and its complement~$\set I^ \complement$ the weight function is given by
\begin{equation} \label{eq:oneDWeightChoice}
	w^B(x) = \sum_{i \in I^\complement} B_i(x) \, .
\end{equation}
The so-assembled weight functions have the following properties.
\newtheorem{theorem}{Proposition}
\begin{theorem}
	The weight functions~$w^B(x)$ and~$w^Q(x)$ are at least~$C^1$ smooth, piecewise polynomials in the parameter space, have local support and form a partition of unity.
\end{theorem}

Finally, using the above set of weight and basis functions, we define the smooth blended B-splines, or SB-splines, as
\begin{equation} \label{eq:oneDBlendedBasisFunctions}
	\vec N(x) = \begin{pmatrix}
		w^B(x) \vec B(x)^\trans & w^Q(x) \vec Q(x)^\trans 
	\end{pmatrix}^\trans 
	= \begin{pmatrix}
		N_1(x) & \dotsc & N_{n_N}(x) 
	\end{pmatrix}^\trans  \, ,
\end{equation}
where~$n_N = n_B + n_Q$. Evidently, the smoothness, the support size and the polynomial degree of the SB-splines depend on the properties of~$\vec B(x)$, $\vec Q(x)$, $w^B(x)$ and $w^Q(x)$. Critical for the smoothness properties of the SB-splines is the choice of the weight function~$w^{B} (x)$. 
\begin{theorem}
	The SB-splines~$\vec N(x)$ are at least~$C^1$ smooth, linearly independent, non-negative, have local support and form a partition of unity.
\end{theorem}
The properties of smoothness, non-negativity, local support and partition of unity follow directly from the blending construction since both $\vec B (x)$ and $\vec Q (x)$ possess these properties.  For the proof of linear independence, see~\ref{sec:independence}.

\begin{figure}[]
	\centering
	\subfloat[][$\vec B(x)$ in black and blue and $w^B(x)$ in red \label{fig:oneDQuadratica}] {
		\includegraphics[width=0.36 \textwidth]{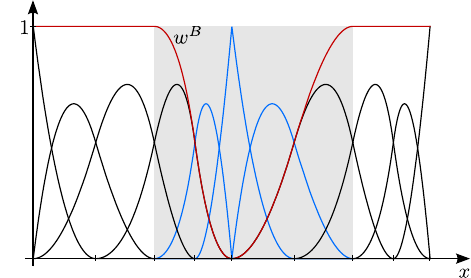} }
	\hfil
		\subfloat[][$\vec Q(x)$ in black and $w^Q(x)$ in red \label{fig:oneDQuadraticb}] {
		\includegraphics[width=0.36 \textwidth]{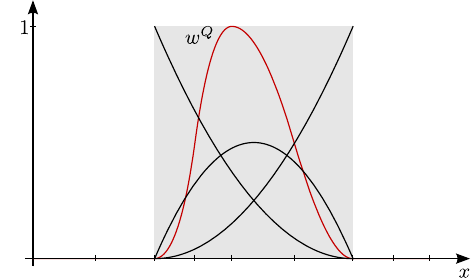} }
	\\
	\subfloat[][$w^B(x) \vec B(x)$  \label{fig:oneDQuadraticc}] {
		\includegraphics[width=0.36 \textwidth]{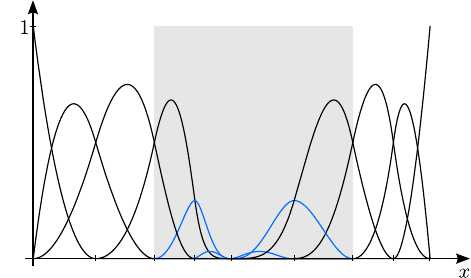} }
	\hfil
	\subfloat[][$w^Q(x) \vec Q(x)$  \label{fig:oneDQuadraticd}] {
		\includegraphics[width=0.36 \textwidth]{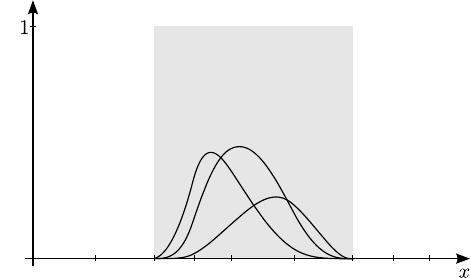} }
	\caption{Blending of the mixed smoothness quadratic B-splines~$\vec B(x)$ in (a) with the quadratic Bernstein basis functions~$\vec Q(x)$ in (b). The non-uniform open knot vector for the B-splines~$\vec B(x)$ has repeated knot values at the boundaries and at the extraordinary point at the centre of the physical domain~$\Omega$. In (a) the B-splines in blue are $C^0$ continuous at the extraordinary point. The blending domain~$\Omega^Q$ is shaded in grey. The weight function~$w^B( x )$ in (a) is the sum of the B-spline basis functions which are $C^1$~smooth at the extraordinary point. Its complement to one is the weight function~$w^Q( x )$ in (b).\label{fig:oneDQuadratic}}
	
	\bigskip
	
	\subfloat[][$\vec B(x)$ in black and blue and $w^B(x)$ in red \label{fig:oneDCubica}] {
		\includegraphics[width=0.36 \textwidth]{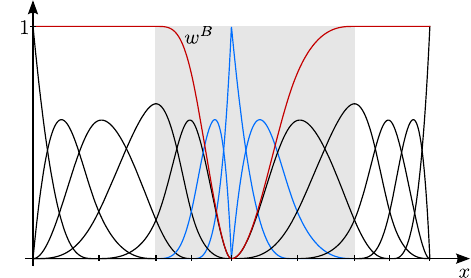} }
	\hfil
	\subfloat[][$\vec Q(x)$ in black and $w^Q(x)$ \label{fig:oneDCubicb}] {
		\includegraphics[width=0.36 \textwidth]{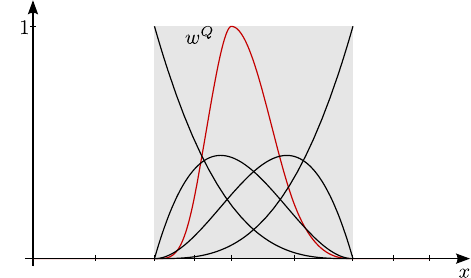} }
	\\
	\subfloat[][$w^B(x) \vec B(x)$ \label{fig:oneDCubicc}] {
		\includegraphics[width=0.36 \textwidth]{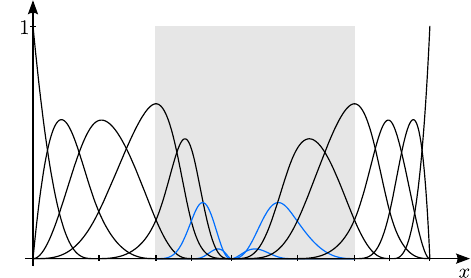} }
	\hfil
	\subfloat[][$w^Q(x) \vec Q(x)$ \label{fig:oneDCubicd}] {
		\includegraphics[width=0.36 \textwidth]{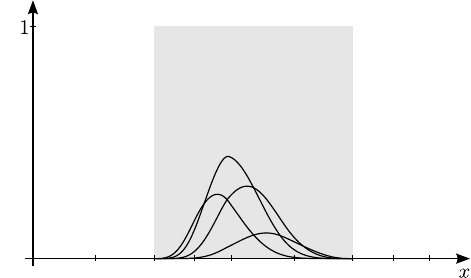} }
	\caption{Blending of the mixed smoothness cubic B-splines~$\vec B(x)$ in (a) with the cubic Bernstein basis functions~$\vec Q(x)$ in (b). The non-uniform open knot vector for the B-splines~$\vec B(x)$ has repeated knot values at the boundaries and at the extraordinary point at the centre of the physical domain~$\Omega$. In (a) the B-splines in blue are $C^0$ continuous at the extraordinary point. The blending domain~$\Omega^Q$ is shaded in grey. The weight function~$w^B( x )$ in (a) is the sum of the B-spline basis functions which are at least $C^1$~smooth at the domain centre. Its complement to one is the weight function~$w^Q( x )$ in (b). \label{fig:oneDCubic}}
	
\end{figure}

As a concrete example, Figure~\ref{fig:oneDQuadratic} illustrates the blending of mixed smoothness quadratic B-splines~$\vec B(x)$ with  quadratic Bernstein basis functions~$\vec Q(x)$. The physical domain~$\Omega$  has at its centre an extraordinary point with~$C^0$ continuity, introduced using a non-uniform open knot vector for~$\vec B(x)$. Except at the extraordinary point the quadratic B-spline basis~$\vec B(x)$ is $C^1$~smooth. The weight function $w^B( x )$ is assembled from  B-splines~$B_i(x)$ by excluding the ones which are  at most $C^0$~smooth at the extraordinary point. Hence, the weight function $w^B( x )$ and, in turn, its complement \mbox{$w^Q( x ) = 1 -w^B( x )$}  are intrinsically $C^1$~smooth. The blending domain~$\Omega^Q$ is equal to the support of the weight function~$w^{Q}(x)$. After determining the weight functions it is straightforward to compute the basis functions~$ N_i(x)$ depicted in Figures~\ref{fig:oneDQuadraticc} and~\ref{fig:oneDQuadraticd}, which are all $C^1$~smooth.  We emphasise that a key aspect of our construction is that the weight function $w^B( x )$ is assembled from the smooth B-splines~$\vec B(x)$.  As apparent in Figure~\ref{fig:oneDQuadratica}, outside the blending region~$\Omega^Q$ the weight function~$w^B( x )$ is equal to one so that the SB-splines are equal to the standard B-splines.

Evidently, the proposed construction can be applied to B-splines of any degree. Figure~\ref{fig:oneDCubic}  illustrates the blending of mixed smoothness cubic  B-splines and cubic Bernstein basis functions. The SB-splines depicted in Figures~\ref{fig:oneDCubicc} and~\ref{fig:oneDCubicd} are also $C^1$ continuous. In the next sections, we restrict our attention to~$p_B = p_Q = 2$.
\section{Two-dimensional quadratic SB-splines \label{sec:two}}
%
We are given an unstructured quadrilateral mesh describing a domain $\Omega \subset \mathbb R^2$. The mesh consists of elements, i.e. quadrilateral faces, and their edges and vertices. We assume that all the vertices on the boundary of the mesh are regular, i.e. are adjacent to two elements, and all the extraordinary vertices within the mesh, i.e. vertices with other than four adjacent elements, are sufficiently separated as to be specified.   If not the case, this can be achieved by successive quadrisection refinement of all elements. The new vertices introduced during refinement are all regular so that the extraordinary vertices become more and more separated.  In the following, without loss of generality, we focus on bi-quadratic B-splines and assume that the mesh has a single extraordinary vertex of valence~$v \neq 4$.  In the (locally) structured regions of the mesh standard smooth tensor-product bi-quadratic B-splines can be defined. It is impossible to define such B-splines in the $1$-neighbourhood of an extraordinary vertex due to the lack of a tensor-product  mesh structure. A $1$-neighbourhood of a vertex is formed by the union of the elements that contain the vertex. The $n$-neighbourhood is defined recursively as the union of all 1-neighbourhoods of all the vertices in the $(n-1)$-neighbourhood. With this definition at hand, we require that the 3-neighbourhoods of the extraordinary vertices in the considered mesh are disjoint.

\subsection{Review of mixed B-splines\label{sec:C0blended}}
%
Although it is not possible to define a standard tensor-product B-spline basis on an unstructured mesh, it is possible to construct a B-spline basis of mixed smoothness; see Toshniwal~\cite{toshniwal2022quadratic}. The mixed B-splines are~$C^1$ continuous away from the~$1$-neighbourhood of extraordinary vertices and are~$C^0$ continuous along mesh edges adjacent to extraordinary vertices. That is, away from the~$1$-neighbourhood of extraordinary vertices the mixed B-splines are identical to tensor-product B-splines. On structured meshes there is a one-to-one correspondence between the bi-quadratic B-splines and elements (away from the boundaries). This is also the case for mixed B-splines. Hence, we can assign a control vertex to each element. The support of a mixed B-spline consists of all elements sharing a vertex with the respective element.

\begin{figure}[]
	\centering
	\subfloat[][Mixed B-spline \& B\'ezier control vertices \label{fig:twoDBezierAveragingIntro}] {
		\includegraphics[width=0.3 \textwidth]{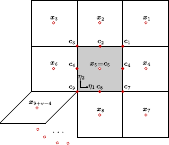} }
	\hfill
	\subfloat[][Edge mask \label{fig:twoDBezierAveragingEdge}] {
		\includegraphics[width=0.3 \textwidth]{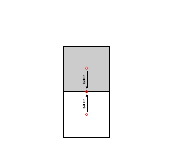} }
	\hfil
	\subfloat[][Vertex mask for valence~$v$\label{fig:twoDBezierAveragingVertex}] {
		\includegraphics[width=0.3 \textwidth]{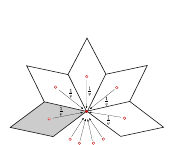} }
	\caption{Averaging masks for computing the bi-quadratic B\'ezier control vertices~$\vec c_j$.  The empty circles denote the mixed B-spline control vertices~$\vec x_i$ and the solid circles the B\'ezier control vertices~$\vec c_j$. The corner mask (c) for extraordinary vertices is a straightforward generalisation of the corner mask for ordinary vertex with~$v = 4$. The masks (b) and (c) describe a bi-quadratic B-spline when an element's all vertices are regular.  \label{fig:twoDBezierAveraging}}
\end{figure}

We represent the non-zero mixed B-splines within an element with bi-quadratic B\'ezier basis functions. The control vertices of the mixed B-splines are denoted with \mbox{$\vec x_i \in \mathbb R^2$} and the ones of the B\'ezier basis functions with~\mbox{$\vec c_j  \in \mathbb R^2$.} The numbering of both sets of control vertices is given in Figure~\ref{fig:twoDBezierAveragingIntro}. We define the mixed B-splines by first establishing the map from the control vertices $\vec x_i$ to~$\vec c_j$. To this end, the B\'ezier control vertices are expressed as linear combinations of  the mixed B-spline control vertices. The corresponding weights can be graphically visualised with the masks shown in Figures~\ref{fig:twoDBezierAveragingEdge} and~\ref{fig:twoDBezierAveragingVertex}.  The edge B\'ezier control vertices $\vec c_2$, $\vec c_4$, $\vec c_6$ and $\vec c_8$ are determined using the mask in Figure~\ref{fig:twoDBezierAveragingEdge} and the corner B\'ezier control vertices $\vec c_1$, $\vec c_2$, $\vec c_3$ and $\vec c_4$ using the mask in Figure~\ref{fig:twoDBezierAveragingVertex}. The centre B\'ezier control vertex~$\vec c_5$ has the same value as the mixed B-spline control vertex~$\vec x_5$. Finally, the mapping of the mixed B-spline control vertices~$\vec x_i$ to the B\'ezier control vertices~$\vec c_j$ is given by
\begin{equation} \label{eq:cTox}
\begin{pmatrix}
	\vec c_1 \\ \vec c_2 \\ \vec c_3 \\ \vec c_4 \\ \vec c_5 \\ \vec c_6 \\ \vec c_7 \\ \vec c_8 \\ \vec c_9
\end{pmatrix}
= 
\begin{pmatrix}
\frac{1}{4} & \frac{1}{4} & 0 & \frac{1}{4} & \frac{1}{4} & 0 & 0 & 0 & \cdots & 0
\\ 0 & \frac{1}{2} & 0 & 0 & \frac{1}{2} & 0 & 0 & 0 & \cdots & 0
\\ 0 & \frac{1}{4} & \frac{1}{4} & 0 & \frac{1}{4} & \frac{1}{4} & 0 & 0 & \cdots & 0
\\ 0 & 0 & 0 & \frac{1}{2} & \frac{1}{2} & 0 & 0   & 0 & \cdots & 0
\\ 0 & 0 & 0 & 0 & 1 & 0 & 0 & 0 & \cdots & 0
\\ 0 & 0 & 0 & 0 & \frac{1}{2} & \frac{1}{2} & 0   & 0 & \cdots & 0
\\ 0 & 0 & 0 & \frac{1}{4} & \frac{1}{4} & 0 & \frac{1}{4} & \frac{1}{4} & \cdots & 0
\\ 0 & 0 & 0 & 0 & \frac{1}{2} & 0 & 0 & \frac{1}{2} & \cdots & 0
\\ 0 & 0 & 0 & 0 & \frac{1}{v} & \frac{1}{v} & 0 & \frac{1}{v} & \cdots & \frac{1}{v}
\end{pmatrix}
\begin{pmatrix}
	\vec x_1 \\ \vec x_2 \\ \vec x_3 \\ \vec x_4 \\ \vec x_5 \\ \vec x_6 \\ \vec x_7 \\ \vec x_8 \\ \vdots  \\ \vec x_{9 + v - 4}
\end{pmatrix} 
\qquad 
	\Rightarrow 
\quad \vec c_j = \sum_i M_{ji} \vec x_i
\, .
\end{equation}

Combining the B\'ezier basis functions~$Q_j(\vec \eta)$ of an element with the map~\eqref{eq:cTox} from the control vertices~$\vec x_i$ to~$\vec c_j$, we can define both a (local) parametrisation of the physical domain~$\Omega$ and the mixed B-splines~$B_i(\vec x)$. Note that the B\'ezier basis functions~$Q_j(\vec \eta)$ are local to each element. Specifically, with the mixed B-spline control vertices~$\vec x_i$ given in Figure~\ref{fig:twoDBezierAveragingIntro} and the map~\eqref{eq:cTox}  the geometry parametrisation within the element corresponding to the control vertex~$\vec x_5$ is given by
\begin{equation} \label{eq:2dParam}
	\vec x(\vec \eta) =  \sum_{j=1}^9 Q_j(\vec \eta) \vec c_j = \sum_{j=1}^9  \sum_{i = 1}^{9+v-4} Q_j(\vec \eta)  M_{ji}  \vec x_i  \quad~\text{with } \vec \eta = (\eta_1, \, \eta_2)  \in \Box:=[0,1] \times [0,1] \,.
\end{equation}
This description also defines the mixed B-spline~$B_5(\vec x)$ associated to control vertex~$\vec x_5$. According to~\eqref{eq:2dParam}, its preimage in the parametric domain is given by 
\begin{equation}
	B_5(\vec \eta) = \sum_{j=1}^9 Q_j (\vec \eta) M_{j5} \, 
\end{equation}	
such that
 \begin{equation} \label{eq:2dMapping}
 	B_5(\vec x) = B_5(\vec \eta) \circ \vec x(\vec \eta)^{-1} \, . 
 \end{equation}
The B-splines~$B_i(\vec x) $ associated to the other control vertices in the mesh are obtained in the same way. See~\cite{toshniwal2022quadratic} for other properties of the mixed B-splines.

\begin{figure}[]
	\centering
	\subfloat[][Quadrilateral mesh \label{fig:twoDQuadMeshInput}] {
		\includegraphics[width=0.3 \textwidth]{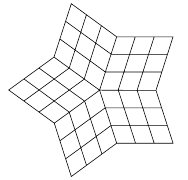} }
	\hfill
	\subfloat[][B\'ezier mesh \label{fig:twoDQuadMeshOutput}] {
		\includegraphics[width=0.3 \textwidth]{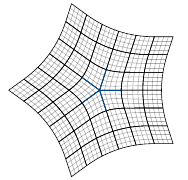} }
	\hfill
	\subfloat[][Basis function \label{fig:twoDQuadMeshSpline}] {
		\includegraphics[width=0.3 \textwidth]{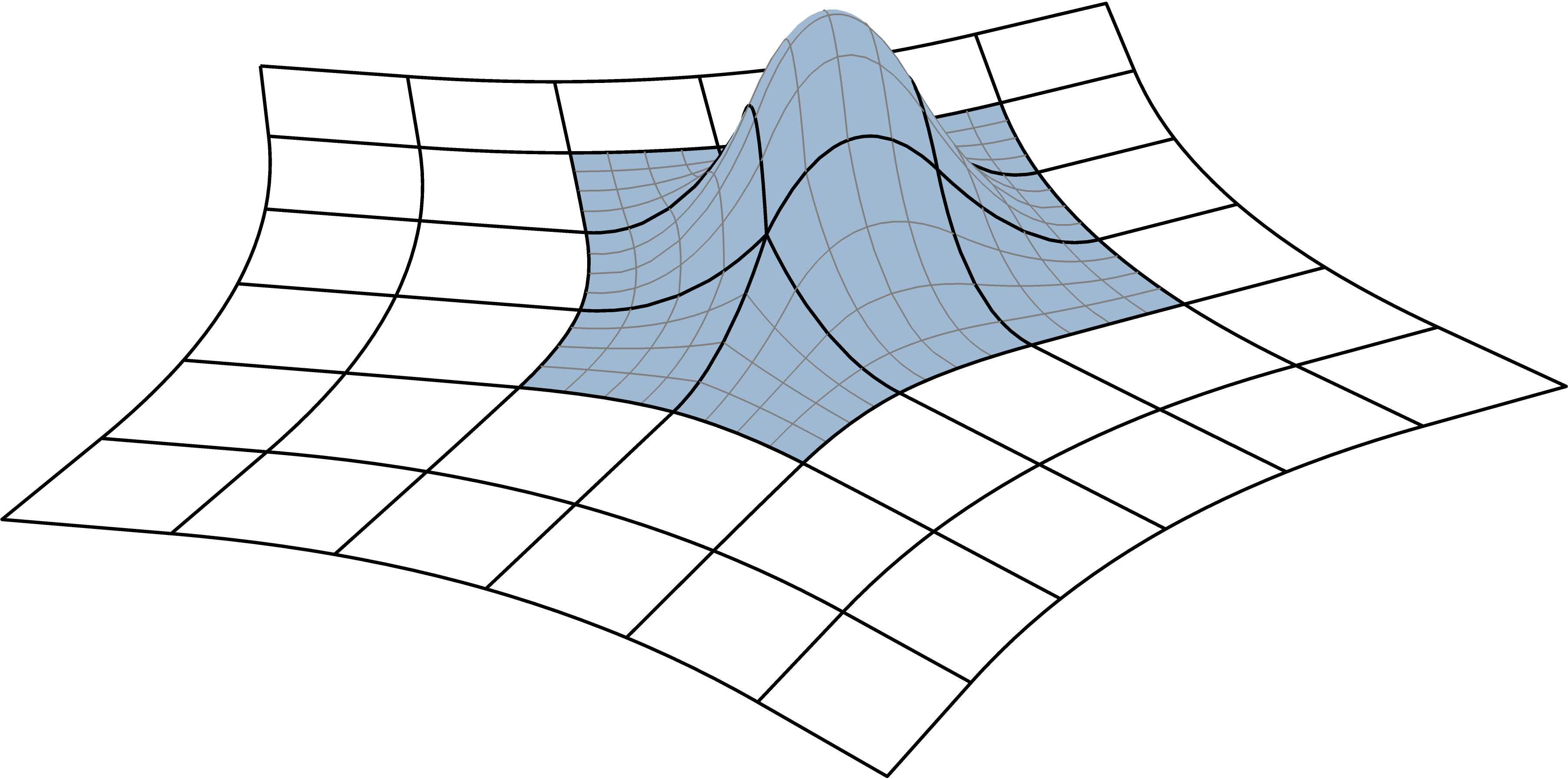} }
	\caption{Unstructured quadrilateral mesh with only one extraordinary vertex with valence~\mbox{$v=5$}, its parametrisation and the graph of one of the corresponding basis functions. The thin lines in (b) and (c) indicate the parameter lines with either \mbox{$\eta_1=\text{const.}$} or \mbox{$\eta_2=\text{const.}$} In the 1-neighbourhood of the extraordinary vertex across the element edges, i.e.\ the blue edges in (b), the parametrisation is only $C^0$ continuous and everywhere else it is $C^1$~smooth.   \label{fig:twoDQuadMesh}}
\end{figure}

As an example, the parametrised domain, i.e.\ B\'ezier mesh corresponding to the unstructured mesh in Figure~\ref{fig:twoDQuadMeshInput} with an extraordinary vertex with~\mbox{$v=5$} is visualised in Figure~\ref{fig:twoDQuadMeshOutput}.  The parametrisation is $C^1$~smooth in most parts of the domain as is suggested by the plotted parameter lines with either \mbox{$\eta_1=\text{const.}$} or \mbox{$\eta_2=\text{const.}$} It is $C^0$ continuous across the edges that contain the extraordinary vertex. For a control vertex in the 1-neighbourhood of the extraordinary vertex we obtain the basis function  shown in Figure~\ref{fig:twoDQuadMeshSpline}. This basis function is only $C^0$ continuous across the edges that contain the extraordinary vertex, as can be inferred from the plotted parameter lines.
%
\subsection{Construction of SB-splines \label{sec:C1blended}}
%
The construction of the blended $C^1$ smooth basis functions is analogous to 1D. First, we choose a weight function~$w^{B} ( \vec x)$ and its complement to one~$ w^Q(\vec x) = 1 - w^{B} ( \vec x)$. The weight function~$w^{B} ( \vec x)$ is assembled from the smooth mixed B-splines defined on the unstructured mesh. Importantly, the blending of the basis functions takes place in the physical domain~$\Omega$ rather than the parametric domain of the basis functions. It is impossible to map every  multi-dimensional physical domain with arbitrary topology onto a single parametric domain. This can be achieved only by introducing an atlas  consisting of several charts with respective parametric domains and transition functions~\cite{majeedCirak:2016,zhang2020manifold,zhang2019manifold}. The definition of such smooth transition function on unstructured meshes is usually very challenging. Instead, constructing the weight functions on the physical domain sidesteps the need for an atlas and smooth transition functions.

\begin{figure}[t]
	\centering
	\subfloat[][B\'ezier mesh\label{fig:twoDWeightBSplineRings}] {
		\includegraphics[width=0.3 \textwidth]{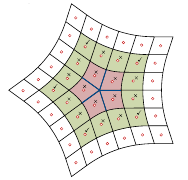} }
	\hfill
	\subfloat[][$w^B(\vec x)$\label{fig:twoDWeightwB}] {
		\includegraphics[width=0.3 \textwidth]{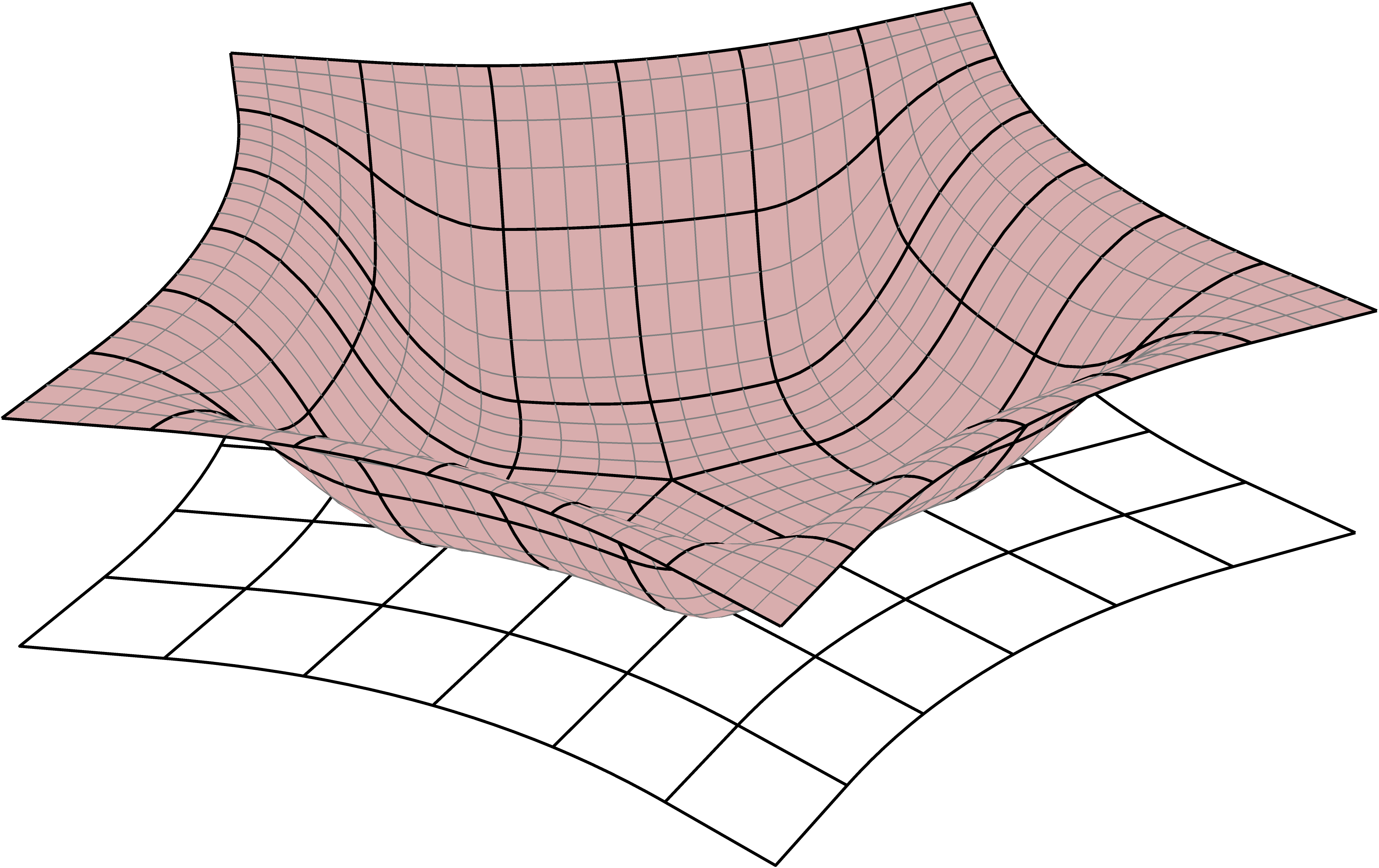} }
	\hfill
	\subfloat[][$w^Q(\vec x)$ \label{fig:twoDWeightwQ}] {
		\includegraphics[width=0.3 \textwidth]{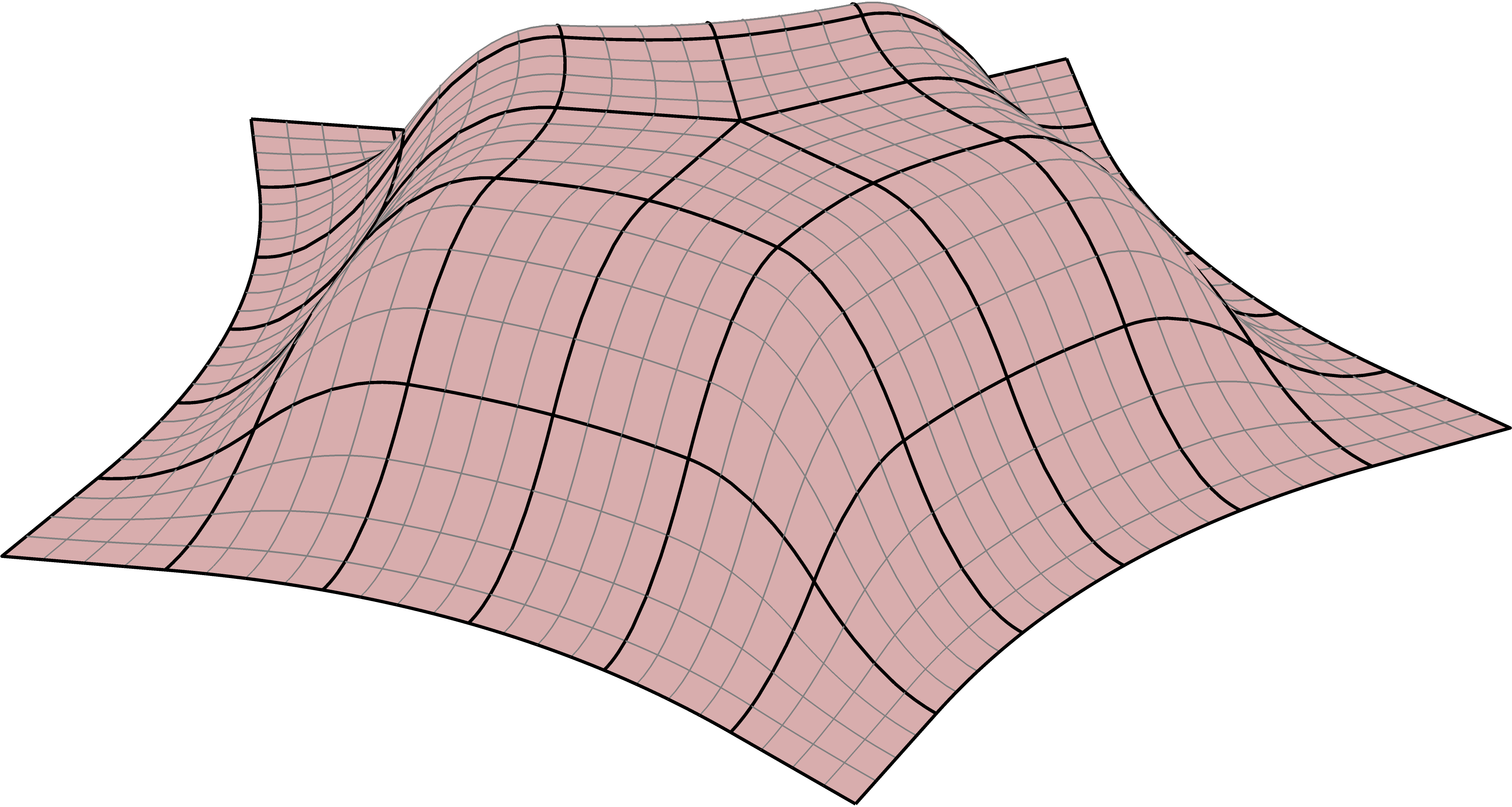} }
	\caption{Smoothness of the mixed B-splines~$B_i(\vec x)$ and weight functions~$w^B(\vec x)$ and~$w^{Q} ( \vec x)$. In (a) the mixed B-splines~$ B_i( \vec x)$ are represented by the respective control vertices (empty circles). In addition,~$ B_i( \vec x)$ that are~$C^1$ and~$C^0$ continuous across element edges in the $1$-neighbourhood of the extraordinary vertex (blue edges) are labelled with ticks and crosses, respectively. In (b) the weight function~$w^B(\vec x)$ is assembled by excluding the mixed B-splines~$B_i(\vec x)$ belonging to the control vertices in the 2-neighbourhood of the extraordinary vertex. Its complement to one~$w^{Q} ( \vec x) = 1 - w^{B} ( \vec x)$ is shown in (c). \label{fig:twoDWeight}}
\end{figure}

As discussed in the preceding section smooth mixed B-splines~$B_i(\vec \eta)$ are defined only away from the \mbox{1-neighbourhood} of an extraordinary vertex. Furthermore, the mixed B-splines~$B_i(\vec x)$ on the physical domain are obtained by mapping~$B_i(\vec \eta)$ from the reference element via the mapping~$\vec x(\vec \eta)^{-1}$, see~\eqref{eq:2dMapping}. According to the chain rule of differentiation the smoothness of~$B_i(\vec x)$ relies both on the smoothness of~$B_i(\vec \eta)$ and~$\vec x(\vec \eta)^{-1}$. This implies that beyond the 1-neighbourhood of an extraordinary vertex  most of the control vertices in its  2-neighbourhood belong to non-smooth mixed B-splines~$B_i(\vec x)$ as well, see Figure~\ref{fig:twoDWeightBSplineRings}.
 
We assemble the weight function~$w^{B} ( \vec x)$ from the mixed B-splines~$B_i(\vec x)$ by excluding the ones belonging to the control vertices in the 2-neighbourhood of the extraordinary vertex, i.e.\ by excluding the mixed B-splines associated to all control vertices marked with a cross or a tick in Figure~\ref{fig:twoDWeightBSplineRings}. Although only the non-smooth mixed B-splines must be excluded, for ease of implementation we exclude some of the smooth mixed B-splines as well. The so-obtained smooth weight function~$w^{B} ( \vec x)$ and its complement to one~$w^{Q} ( \vec x) = 1-w^{B} ( \vec x)$ are depicted in Figures~\ref{fig:twoDWeightwB} and~\ref{fig:twoDWeightwQ}, respectively.   Evidently, both weight functions are $C^1$ smooth, bi-quadratic on the reference element domain~$\Box$   and $\supp w^{Q} ( \vec x)$  is comprised of the 3-neighbourhood of the extraordinary vertex.
\begin{figure}[t]
	\centering
	\subfloat[][$w^B(\vec x) B_i(\vec x)$  \label{fig:twoDBlendedBasisFunsA}] {
		\includegraphics[width=0.3 \textwidth]{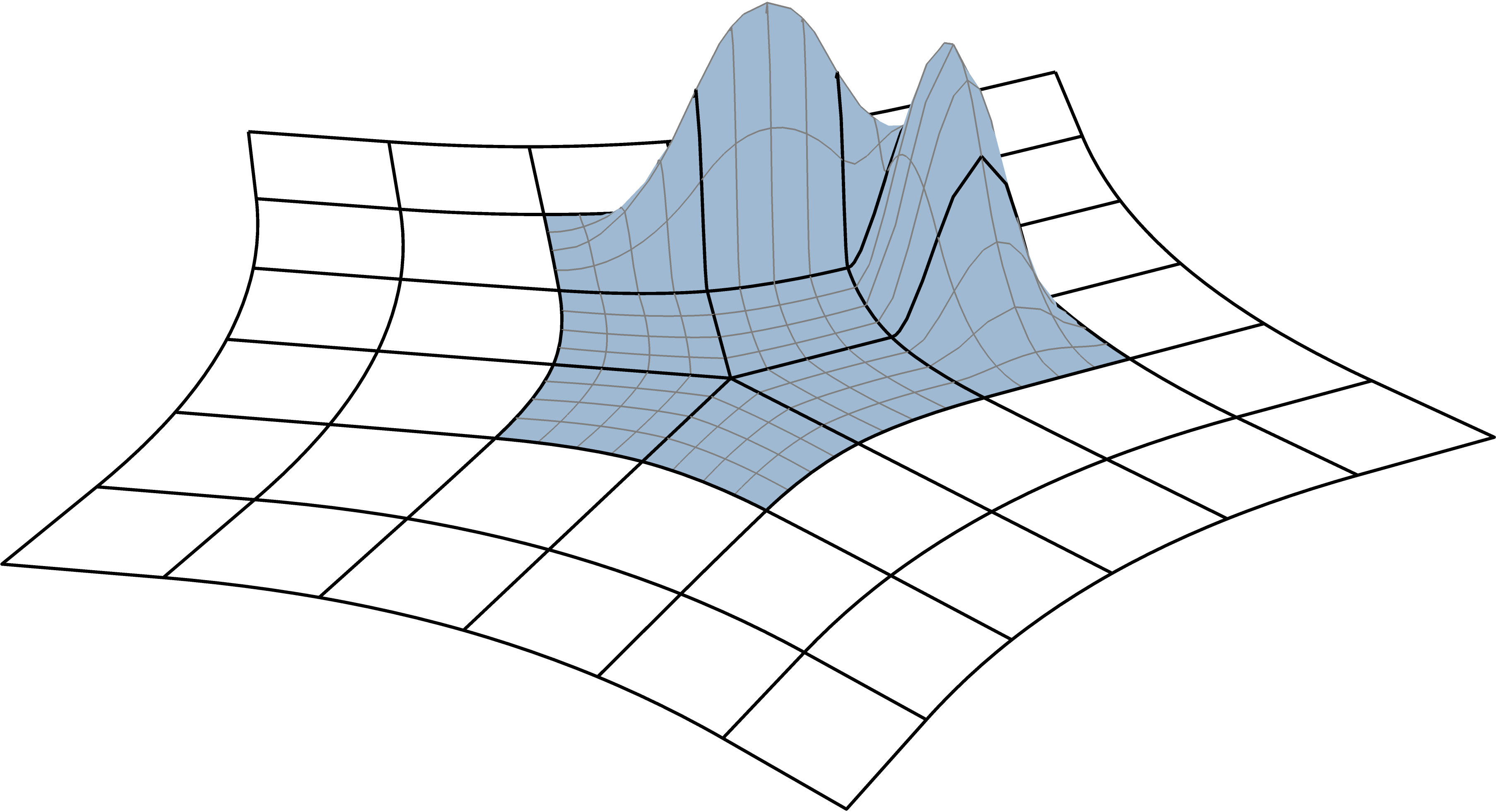} }
	\hfill
	\subfloat[][$w^B(\vec x) B_i(\vec x)$  \label{fig:twoDBlendedBasisFunsB}] {
		\includegraphics[width=0.3 \textwidth]{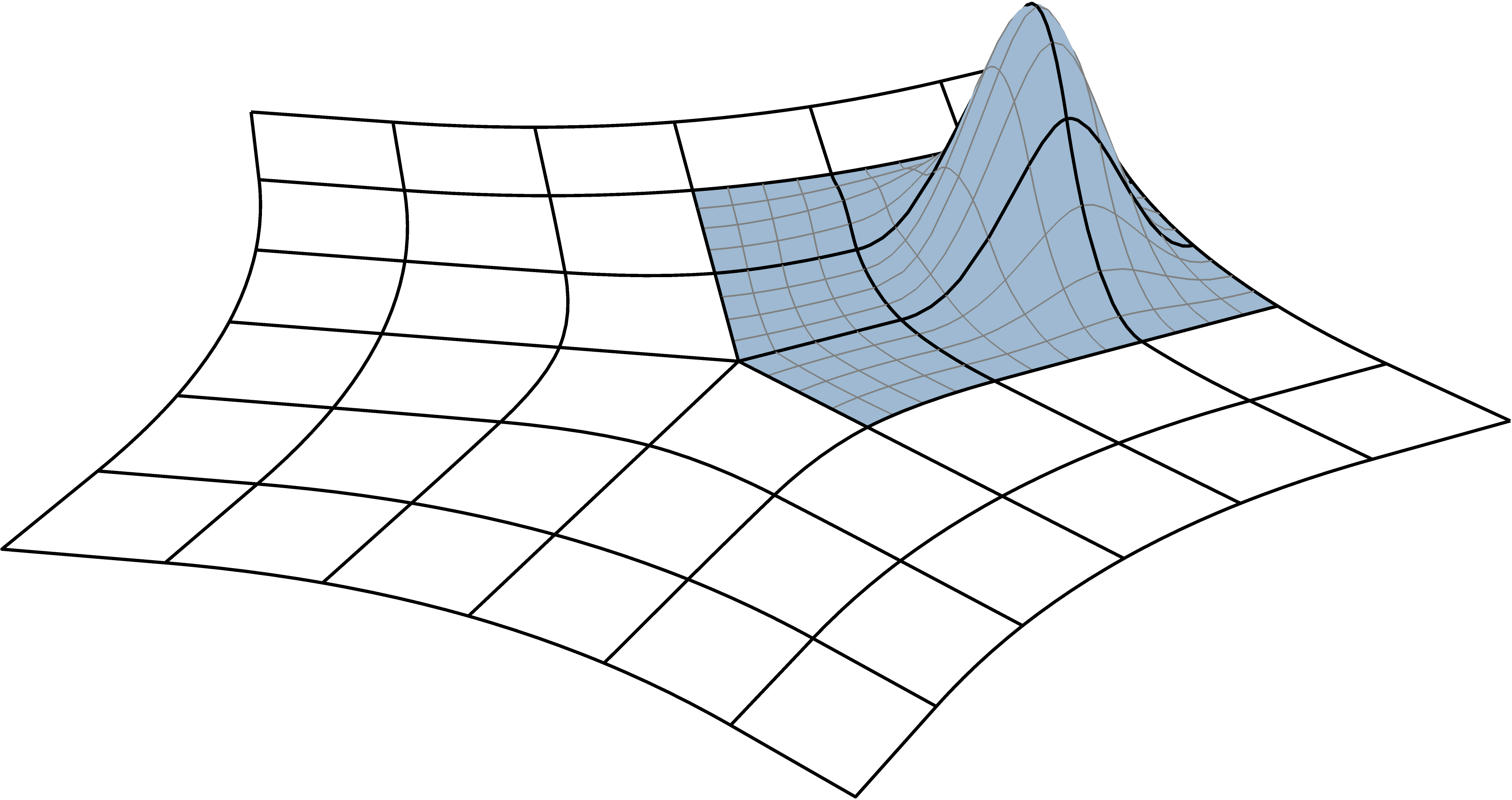} }
	\hfill
	\subfloat[][$w^Q(\vec x)Q_j(\vec x)$    \label{fig:twoDBlendedBasisFunsC}] {
		\includegraphics[width=0.3 \textwidth]{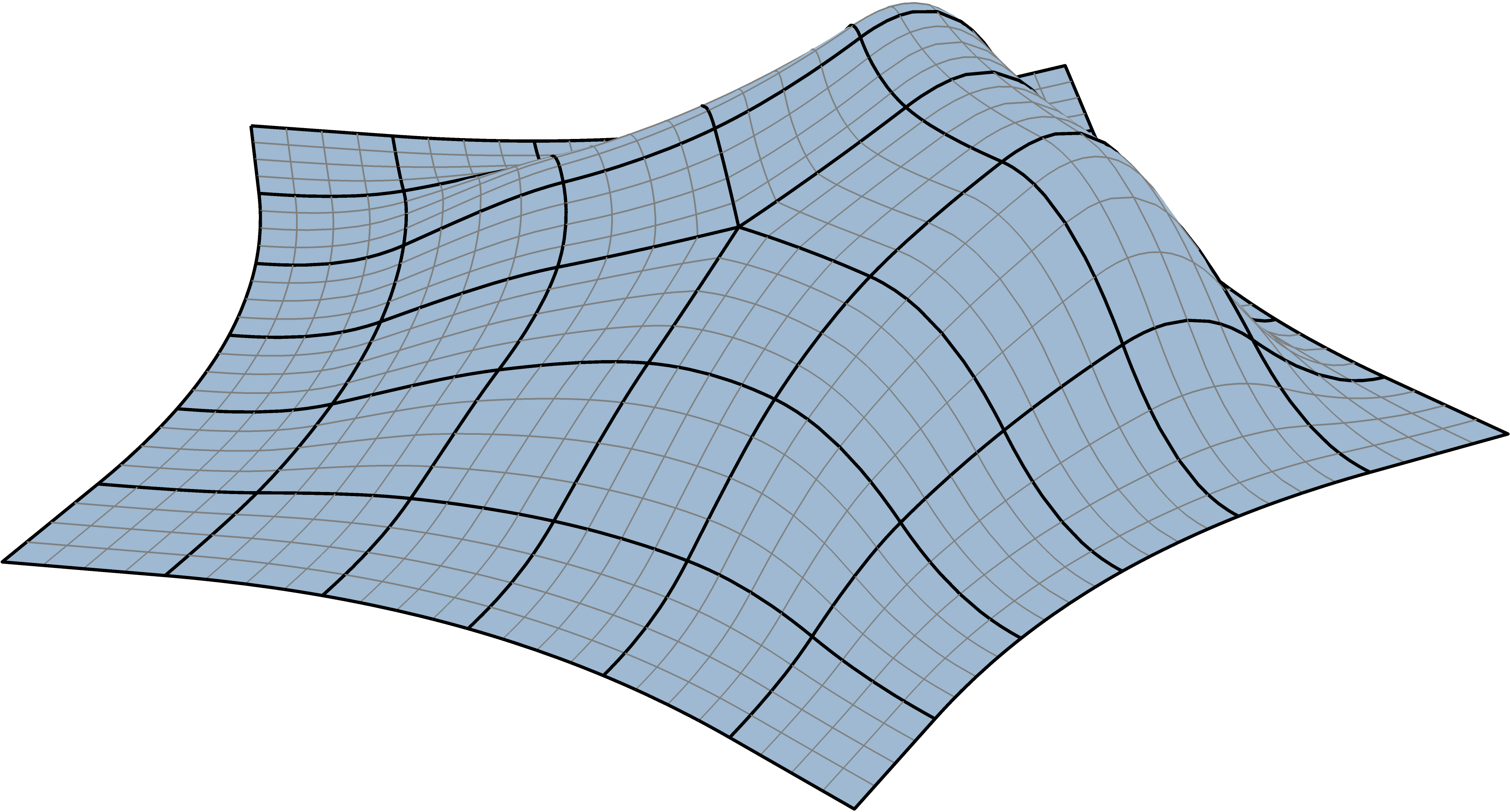} }
	\caption{Three of the obtained SB-splines~$N_i(\vec x)$. The SB-splines in (a) and (b) correspond to the non-smooth mixed B-splines~$B_i(\vec x)$ belonging to the control vertices in the $1$-neighbourhood and $2$-neighbourhood, respectively. The SB-spline in (c) corresponds to one of the B\'ezier basis functions~$Q_j(\vec x)$. \label{fig:twoDBlendedBasisFuns}}
\end{figure}

Next, we choose a bi-quadratic Bernstein basis~$\{ Q_j(\vec x) \}_{j=1}^{9} $ as the second basis for blending. These are defined on a different domain than the ones in~\eqref{eq:2dParam}; indeed, considering that the blending takes place in the physical domain~$\Omega$ this basis is defined in the physical domain with~$\vec x \in \Omega$. To guarantee the positivity of the SB-splines, the domain of the Bernstein basis~$\{ Q_j(\vec x) \}_{j=1}^{9} $ must enclose the blending domain $\Omega^Q$, but can be chosen freely otherwise. As in 1D, the SB-splines are then defined by
\begin{equation}
	\vec N (\vec x)  =
	\begin{pmatrix}
		w^B(\vec x)  \vec B(\vec x)^\trans &
		w^Q(\vec x)  \vec Q (\vec x)^\trans
	\end{pmatrix} ^\trans 
\end{equation} 
with
\begin{equation}
	\vec B (\vec x ) =  \begin{pmatrix} B_1(\vec x) & \dotsc & B_{n_B}(\vec x)  \end{pmatrix}^\trans \, , 
	\quad  \vec Q (\vec x ) = \begin{pmatrix}Q_1(\vec x) & \dotsc & Q_9(\vec x) \end{pmatrix}^\trans \, .
\end{equation}
The obtained basis functions~$\vec N(\vec x)$ are $C^1$~smooth.  In Figures~\ref{fig:twoDBlendedBasisFunsA} and~\ref{fig:twoDBlendedBasisFunsB} two of the basis functions~$w^B(\vec x) B_i(\vec x) $ and in Figure~\ref{fig:twoDBlendedBasisFunsC} one of the basis functions~$w^Q(\vec x) Q_j(\vec x)$  are plotted.

\begin{figure}[]
	\centering
	\subfloat[][Univariate weight functions in red and univariate B-splines in black and blue \label{fig:twoDSectorialWeightUnivariate}] {
		\includegraphics[width=0.4 \textwidth]{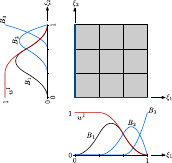} }
	\hfill
	\subfloat[][$w^Q(\vec x)$ for a given sector \label{fig:twoDSectorialWeightSector}] {
		\includegraphics[width=0.4 \textwidth]{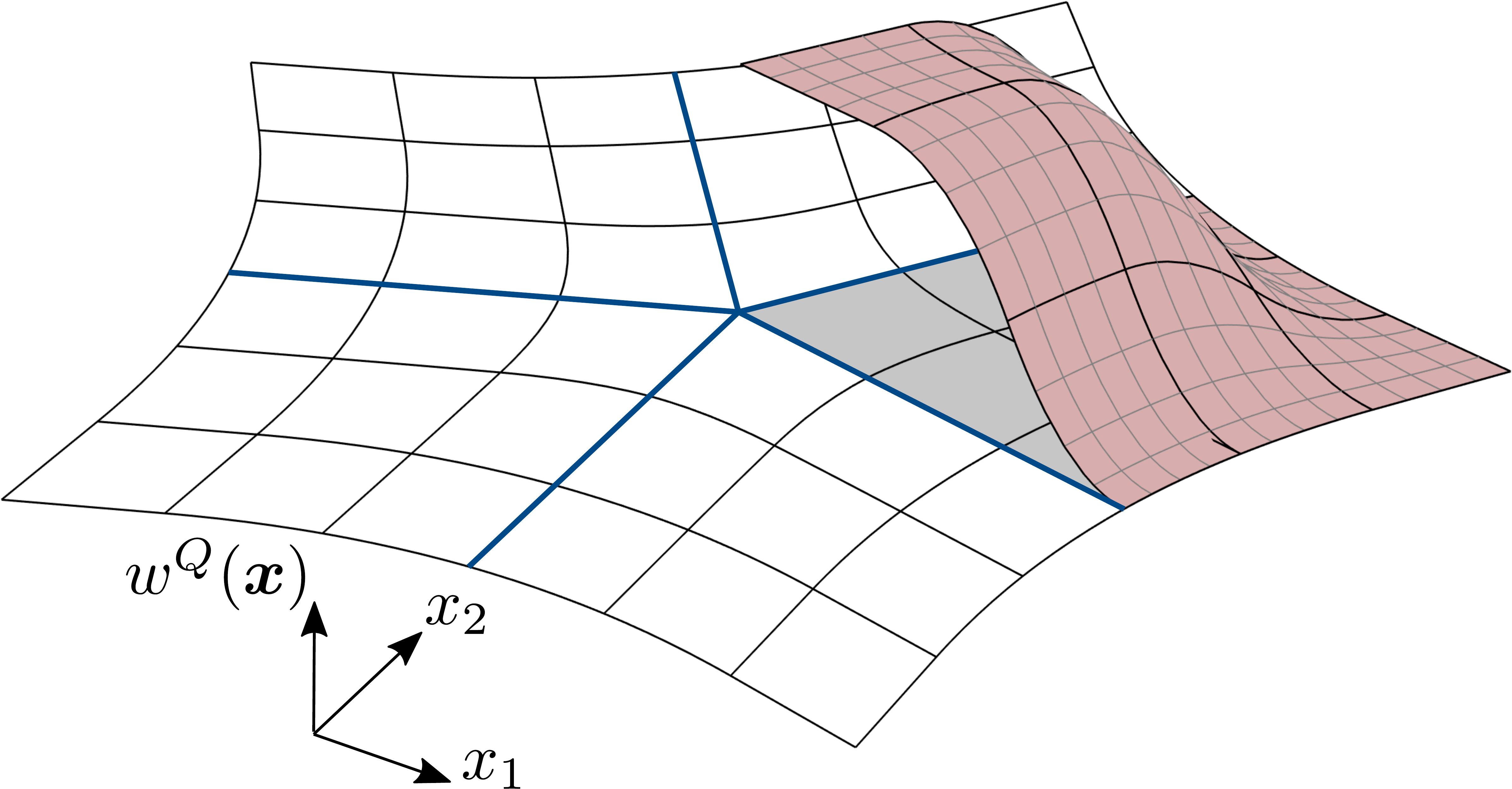} }
	\caption{Construction of the weight function~$w^Q(\vec x)$ on one of the~$v$ sectors in the 3-neighbourhood of an extraordinary vertex. The weight function~$w^Q(\vec \xi)$ on the parametric domain~$\hat \Omega$ is defined as the product of the univariate weight functions~$w^{\rm{I}}(\xi_1)$ and~$w^{\rm{I}}(\xi_2)$ in (a). Formally,~$w^Q(\vec x)$ on each element on the domain~$\Omega$ in (b) is given by \mbox{$w^Q(\vec x) = \vec w^Q(\vec \xi) \circ \vec \eta (\vec \xi)^{-1} \circ \vec x(\vec \eta)^{-1}$}.   \label{fig:twoDSectorialWeight}}
\end{figure}

In usual finite element implementations integrals are evaluated in a reference element domain~\mbox{$\Box:=[0,1] \times [0,1]$}. To facilitate the element-based  implementation of the proposed blended approach, we consider a sector-wise construction of the weight functions~$w^B(\vec x)$ and~$w^Q(\vec x)$. The process is outlined in Figure~\ref{fig:twoDSectorialWeight}. The $3$-neighbourhood of the extraordinary vertex is partitioned into~$v$ sectors. Each sector consisting of~$3 \times 3$ elements is parametrised with~\mbox{$\vec \xi = (\xi_1, \xi_2) \in \hat\Omega:=[0,1] \times [0,1]$}, see Figure~\ref{fig:twoDSectorialWeightUnivariate}. The elements on the parametric domain~$\hat{\Omega}$ are mapped using a mapping~$\vec \eta (\vec \xi) $ composed of a translation and a scaling to the reference element~$\Box$ for integration. 

With the extraordinary vertex located at the origin~$\vec \xi = (0,0)$ of~$\hat\Omega$, the implementation of the weight functions is identical on all the~$\nu$ sectors. Therefore, it is sufficient to detail the implementation only one of the sectors.  As illustrated in Figure~\ref{fig:twoDSectorialWeightSector},  the weight function~$w^Q(\vec x)$ is constructed by first defining~$w^Q(\vec \xi)$ on the parametric domain~$\hat \Omega$. We define $w^Q(\vec \xi)$ as the tensor product of univariate weight functions~$ w^{\rm{I}}(\xi_1)$ and~$ w^{\rm{I}}(\xi_2)$,
\begin{equation} \label{eq:twoDSectorialWeight}
	w^Q(\vec \xi) = w^{\rm{I}}(\xi_1)  \otimes w^{\rm{I}}(\xi_2)  \,.
\end{equation}
First, we assemble along the $\xi_1$ and $\xi_2$ axes the univariate weight functions for~$w^B(\xi_1)$ and~$w^B(\xi_2)$  by excluding the univariate B-splines which are not $C^1$~smooth at the origin~$\vec \xi = (0, \, 0)$, see Figure~\ref{fig:twoDSectorialWeightUnivariate}. Then, their complements to one yield
\begin{equation} \label{eq:twoDUnivariateWeight}
	w^{\rm{I}}(\xi_i) = 1 - B_2(\xi_i) - B_3(\xi_i)   \, ,  \quad i = 1, \, 2 \,. 
\end{equation}

For evaluating the finite element integrals the weight function values at the quadrature points in the reference element~$\Box$ are required, which are obtained from
\begin{equation}
	w^Q(\vec \eta) =  w^Q(\vec \xi) \circ \vec \eta (\vec \xi)^{-1} \, . 
\end{equation}
Here, the mapping~$ \eta (\vec \xi) $ is, as mentioned above, composed of a translation and scaling and can be easily inverted.

\section{Three-dimensional quadratic SB-splines \label{sec:three}}
%
The proposed construction of SB-splines can also be extended to unstructured hexahedral meshes describing a domain~$\Omega \subset \mathbb{R}^3$. A mesh consists of elements, i.e. hexahedral cells, their quadrilateral faces, edges and vertices. In 3D there are in addition to extraordinary vertices also extraordinary edges,  see~Figure~\ref{fig:introSphereExtraordinary}. Extraordinary edges are connected to two extraordinary vertices and regular edges to two regular vertices. In this paper, we assume that there are no edges with one ordinary and one extraordinary vertex and that all edges on the boundary of the domain are regular, i.e. are adjacent to two elements. Furthermore, the valence of an edge~$e$ is defined as the number of elements that share the same two vertices like the edge.

\begin{figure}[]
	\centering
	\subfloat[][3-neighbourhood of all the extraordinary vertices\label{fig:threeDSphereExtraordinary}] {
		\includegraphics[width=0.3 \textwidth]{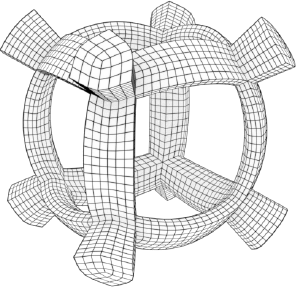} }
	\hfill
	\subfloat[][Extraordinary prisms \label{fig:threeDSpherePrism}] {
		\includegraphics[width=0.3 \textwidth]{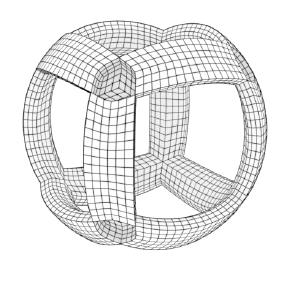} }
	\hfill
	\subfloat[][Extraordinary joints \label{fig:threeDSphereJoint}] {
		\includegraphics[width=0.3 \textwidth]{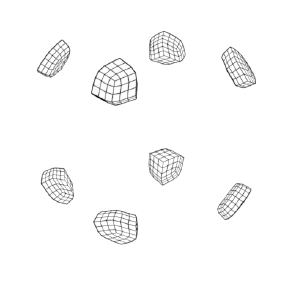} }
	\caption{Extraordinary features in the unstructured hexahedral sphere mesh in Figure~\ref{fig:introSphere}. The 3-neighbourhood of the extraordinary vertices  (a) are decomposed into twenty extraordinary prisms (b) and eight extraordinary joints (c). For illustration purposes, only twelve of the twenty extraordinary prisms are depicted in (b). The extraordinary vertices corresponding to the joints have~$v=4$ and the extraordinary edges corresponding to the prisms have~$e=3$.\label{fig:threeDSphere}}
\end{figure}

Well-designed hexahedral finite element meshes consist of a small number of chains of extraordinary edges. There are usually only a few extraordinary vertices with more than two attached extraordinary edges~\cite{nieser2011cubecover,zhang2018geometric}. As in 2D, only the 3-neighbourhood of the extraordinary vertices and extraordinary edges is relevant for the proposed construction.  In 3D, the union of the 3-neighbourhoods of all the extraordinary vertices in the mesh form a 6-element wide chain of elements as depicted in Figure~\ref{fig:threeDSphereExtraordinary}. We split the chain into several disjoint sets and refer to them as extraordinary prisms or extraordinary joints as illustrated in Figures~\ref{fig:threeDSpherePrism} and \ref{fig:threeDSphereJoint}. Joints consist of the 3-neighbourhood of extraordinary vertices where more than two extraordinary edges meet. The remaining elements in the chain form the prisms. Each prism is connected to either a joint or the domain boundary.

In practice, the possible number of extraordinary edges meeting at an extraordinary vertex is limited. For the sake of clarity and conciseness, without loss of generality, we consider in this section only a joint with four attached prisms, i.e.~$v=4$, each of which have valence~$e=3$. The arbitrary~$v$ and~$e$ case can be similarly elaborated upon as is briefly discussed in~\ref{sec:joint}.

%
\subsection{Review of mixed B-splines \label{sec:C0blendedHex}}
%
We now outline the extension of the mixed B-spline construction in Section~\ref{sec:C0blended} to an unstructured hexahedral mesh. The resulting tri-quadratic mixed B-splines are~$C^1$ continuous away from the~$1$-neighbourhood of extraordinary edges and are only~$C^0$ continuous along mesh faces adjacent to extraordinary vertices.  Again, there is a one-to-one correspondence between the tri-quadratic mixed B-splines and the elements in the mesh (away from the boundaries) so that we assign a control vertex to each element. The support of the respective mixed B-spline consists of all elements sharing a vertex with the element.

\begin{figure}[]
	\centering
	\subfloat[][Face mask \label{fig:threeDBezierAveragingA}] {
		\includegraphics[width=0.3 \textwidth]{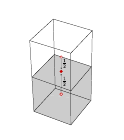} } 
	\hfil
	\subfloat[][Edge mask for edge valence~$e = 4$ \label{fig:threeDBezierAveragingB}] {
		\includegraphics[width=0.3 \textwidth]{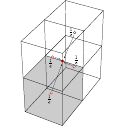} }
	\hfil
	\subfloat[][Corner mask for vertex valence~$v = 8$\label{fig:threeDBezierAveragingC}] {
		\includegraphics[width=0.3 \textwidth]{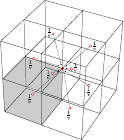} }
	\caption{Averaging masks for computing the tri-quadratic B\'ezier control vertices.  The empty circles denote the mixed B-spline control vertices~$\vec x_i$ and the solid circles the B\'ezier control vertices~$\vec c_j$. The averaging masks for the edge B\'ezier control vertex in (b) and corner B\'ezier control vertex in (c) depend on the valence of the edge~$e$ and the valence of the vertex~$v$, respectively. \label{fig:threeDBezierAveraging}}
\end{figure}

As before, we first represent the non-zero mixed B-splines within an element with tri-quadratic B\'ezier basis functions. Subsequently, the corresponding B\'ezier control vertices are expressed as linear combinations of the mixed B-spline control vertices using the masks depicted in Figure~\ref{fig:threeDBezierAveraging}, the labelling of the vertices has been omitted for simplicity. The face B\'ezier control vertices are determined using the mask in Figure~\ref{fig:threeDBezierAveragingA}, the edge B\'ezier control vertices using the mask in Figure~\ref{fig:threeDBezierAveragingB} and the corner B\'ezier control vertices using the mask in~\ref{fig:threeDBezierAveragingC}. The masks for edge and corner B\'ezier control vertices depend on the valence of the edge~$e$ and the valence of the vertex~$v$, respectively. The centre B\'ezier control vertex has the same value as the element's respective mixed B-spline control vertex.

As in the two-dimensional case, combining the obtained B\'ezier basis control vertices with B\'ezier basis function we can define a (local) parametrisation of the physical domain~$\Omega$ as well as mixed B-splines. Recall that the parameterisation for the quadrilateral mesh in Section~\ref{sec:C0blended} was~$C^0$ continuous across all edges adjacent to the extraordinary vertex. For hexahedral meshes the parametrisation is~$C^0$ continuous across all faces adjacent to the extraordinary edge.
%
\subsection{Construction of SB-splines \label{sec:C1blendedHex}}
%
The construction of the blended~$C^1$~smooth basis functions on unstructured hexahedral meshes follows the 1D and 2D constructions with only slight modification. The key idea is again to consider all the smooth mixed B-splines~$B_i (\vec x)$ to define the weight function~$w^B(\vec x)$ and its complement to  one~\mbox{$w^Q(\vec x) = 1 - w^B(\vec x)$}. A naive implementation of this idea leads on hexahedral meshes to a single weight function~$w^Q (\vec x)$ with a support covering all the connected extraordinary prisms and joints in the mesh. Clearly, such a construction will lead to an overly dense stiffness matrix and adversely affect the approximation properties of the SB-splines. Therefore, as will be detailed in the following, the weight function~$w^Q(\vec x)$  is partitioned into two sets of locally supported weight functions~$w^P_{k, \ell}(\vec x)$  and~$w^J_\ell(\vec x)$ such that
\begin{equation} \label{eq:threeDwQ}
	w^Q(\vec x) =  \sum_{k=1}^{n_P} \sum_{\ell=1}^{n_k} w^{P}_{k,\ell}(\vec x) + \sum_{\ell=1}^{n_J} w^{J}_\ell(\vec x) \, ,
\end{equation}
where~$n_P$ is the number of extraordinary prisms,~$n_k$ is the number of weight functions defined along the corresponding extraordinary prisms and $n_J$ is the number of extraordinary joints. In other words, there is one weight function for each extraordinary joint and several weight functions for each extraordinary prism.  The construction of the prism weight functions~$w^{P}_{k,\ell}(\vec x) $ are introduced in~Section~\ref{sec:weightPrism} and the joint weight functions~$w^{J}_\ell(\vec x)$ in Section~\ref{sec:weightJoint}.

After the weight functions are determined, for blending we assign a tri-quadratic Bernstein basis~$\vec Q(\vec x)$ to  each weight function. In the following, the domain of each tri-quadratic Bernstein basis~$\vec Q(\vec x)$ is assumed to be a cuboid enclosing the support of the corresponding weight function it is assigned to. We index the Bernstein basis similarly to the associated weight functions. Hence, similar to 1D and 2D, the SB-splines are defined by
\begin{equation}
\vec N (\vec x)  =
\begin{pmatrix}
w^B(\vec x) 
\vec B (\vec x )^\trans
&
w^P_{1,1}(\vec x) \vec Q_{1,1}(\vec x)^\trans
&
\cdots
&
w^P_{n_P, n_{n_P}}(\vec x) \vec Q_{n_P, n_{n_P}}(\vec x)^\trans
&
w^J_{1}(\vec x) \vec Q_{1}(\vec x)^\trans
&
\cdots
&
w^J_{n_J}(\vec x) \vec Q_{n_J}(\vec x)^\trans
\end{pmatrix} ^\trans \, .
\end{equation}
%

%
\subsubsection{Weight functions for extraordinary prisms \label{sec:weightPrism}}
\begin{figure}[]
	\centering
	\subfloat[][Extraordinary prism with only two of the sectors shown \label{fig:threeDPrismSegmentA}] {
		\includegraphics[width=0.3 \textwidth]{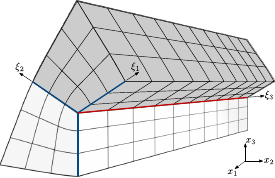} } 
	\hfill
	\subfloat[][Univariate weight functions in red and univariate B-splines in black and blue \label{fig:threeDPrismSegmentB}] {
		\includegraphics[width=0.6\textwidth]{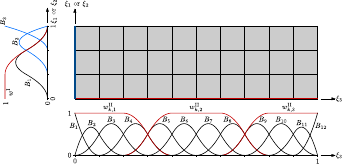} }
	\caption{Construction of the weight function for an extraordinary prism. The two ends of the prism are either joined to an extraordinary joint or the boundary of the domain~$\Omega$.  In (a) only two of the~$e$ sectors are shown for visualisation purposes. Each sector consists of~$3\times 3\times n_{ee}$ elements and has a corresponding (cuboidal) parametric domain~$\hat \Omega$ as illustrated in (b).  The weight functions~$w_{k,\ell}^P(\vec \xi)$ for one of the~$e$ sectors is defined as the tensor product of the bivariate weight function~$w^{\rm{I}}(\xi_1)\otimes w^{\rm{II}}(\xi_2)$ and the three univariate weight functions $w_{k,1}^{\rm{II}}(\xi_3)$, $w_{k,2}^{\rm{II}}(\xi_3)$ and $w_{k,3}^{\rm{II}}(\xi_3)$.  \label{fig:threeDPrismSegment}}
\end{figure}

\begin{figure}[]
	\centering
	\subfloat[][$w^P_{k,1}(\vec x)$ \label{fig:threeDPrismWeightA}] {
		\includegraphics[width=0.2 \textwidth]{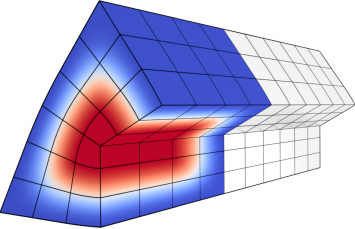} } 
	\hfill
	\subfloat[][$w^P_{k,2}((\vec x)$ \label{fig:threeDPrismWeightB}] {
		\includegraphics[width=0.2 \textwidth]{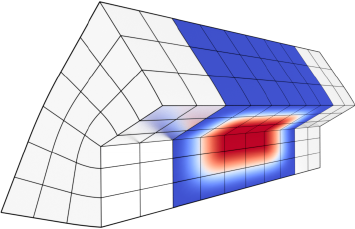} }
	\hfill
	\subfloat[][$w^P_{k,3}((\vec x)$ \label{fig:threeDPrismWeightC}] {
		\includegraphics[width=0.2 \textwidth]{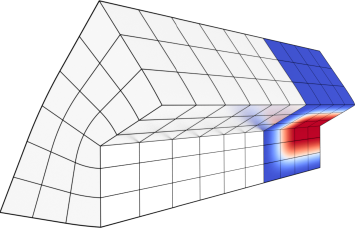} }
	\hfill
	\subfloat[][$w^B(\vec x)$ \label{fig:threeDPrismWeightD}] {
		\includegraphics[width=0.2 \textwidth]{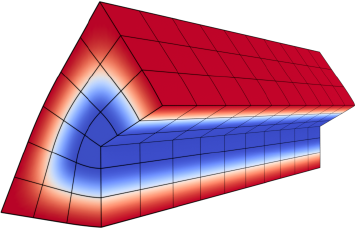} }
	\caption{Weight functions for the extraordinary prism. The weight functions~$w_{k,1}^P(\vec x)$, $w_{k,2}^P(\vec x)$ and $w_{k,3}^P(\vec x)$ are obtained by choosing three univariate weight functions~$w_{k,1}^{\rm{II}} (\xi_3)$, $w_{k,2}^{\rm{II}}(\xi_3)$ and $w_{k,3}^{\rm{II}}(\xi_3)$ along the centre of the extraordinary prism as depicted in Figure~\ref{fig:threeDPrismSegmentB}.  The weight function~$w^B(\vec x)$ is their complement to one, i.e.~$w^B(\vec x) = 1 - w_{k,1}^P(\vec x) - w_{k,2}^P(\vec x) - w_{k,3}^P(\vec x)$. The scalar field ranges between~$0$ (blue) and~$1$ (red). \label{fig:threeDPrismWeight}}
\end{figure}

%
The weight functions for one extraordinary prism are obtained as illustrated in Figures~\ref{fig:threeDPrismSegment} and~\ref{fig:threeDPrismWeight}. The two ends of the prism are either  joined to an extraordinary joint or the boundary of the domain~$\Omega$. The centre of the prism consists of~$n_{ee}$ extraordinary edges of valence~$e=3$. For constructing the weight functions the prism is partitioned into~$e$ sectors, see Figure~\ref{fig:threeDPrismSegmentA}.  Each sector consists of~$3 \times 3 \times n_{ee}$ elements and is parametrised using ~\mbox{$\vec \xi = (\xi_1,\xi_2,\xi_3) \in \hat{\Omega}$} with the parametric domain \mbox{$\hat{\Omega}:=[0,1] \times [0,1] \times [0,1]$}.  The extraordinary edges are located along the parametric axis~\mbox{$\vec \xi = (0, \, 0, \, \xi_3 )$}. The elements on~$\hat{\Omega}$ are mapped to the reference element~$\Box$ for integration using a mapping~$\vec \eta (\vec \xi) $ composed of a translation, rotation and scaling. 

On a given sector of the $k$-th prism, we define~$n_k$ univariate weight functions~$w_{k, \ell}^{\rm{II}}(\xi_3)$ from the available quadratic univariate B-splines~$B_i(\xi_3)$. The number of weight functions~$n_k$ can be chosen flexibly, as long as
\begin{itemize}
	\item[--] each~$w_{k, \ell}^{\rm{II}}(\xi_3)$ is the sum of a certain number of consecutive B-splines,
	\item[--] each~$w_{k,\ell}^{\rm{II}}(\xi_3)$ has vanishing derivatives at the endpoints of its support, 
	\item[--] each B-spline is used to build exactly one~$w_{k, \ell}^{\rm{II}}(\xi_3)$, c.f.\ Figure~\ref{fig:threeDPrismSegmentB},
	\item[--] and the sum of all~$w_{k, \ell}^{\rm{II}}(\xi_3)$ is equal to 1.
\end{itemize}
For instance, in Figure~\ref{fig:threeDPrismSegmentB} the univariate weight functions~$w_{k,\ell}^{\rm{II}}(\xi_3)$ are defined as
\begin{equation}
w_{k,1}^{\rm{II}}(\xi_3) = \sum_{i=1}^4 B_i(\xi_3)\,, \quad  w_{k,2}^{\rm{II}}(\xi_3) = \sum_{i=5}^8 B_i(\xi_3)\,, \quad  w_{k,3}^{\rm{II}}(\xi_3) = \sum_{i=9}^{12} B_i(\xi_3)   \,.
\end{equation}
The isocontours of the three weight functions and their complement to one~$w^B(\vec x)$ on two of the three sectors are depicted in Figure~\ref{fig:threeDPrismWeight}. Note that choosing a large~$n_k$ ensures that the SB-splines have small support sizes.

Following the above, the trivariate weight functions~$w^P_{k,\ell}(\xi_1, \xi_2, \xi_3 )$ are defined as the tensor product of the bivariate weight function~$w^{\rm{I}}(\xi_1) \otimes w^{\rm{I}}(\xi_2)$, introduced in Section~\ref{sec:C1blended}, and the~$n_k$ univariate weight functions~$w_{k, \ell}^{\rm{II}}(\xi_3)$,
\begin{equation} \label{eq:threeDPrismWeight}
	w^P_{k,\ell}(\vec \xi) = w^{\rm{I}}(\xi_1) \otimes w^{\rm{I}}(\xi_2) \otimes w_{k, \ell}^{\rm{II}}(\xi_3) \,, \quad \ell=1, \dotsc, n_k \, .
\end{equation} 
The construction is repeated for all prisms to obtain weight functions for all~$k$.

%
\subsubsection{Weight functions for extraordinary joints \label{sec:weightJoint}}
\begin{figure}[t]
	\centering
	\subfloat[][Extraordinary joint and the attached prisms with some elements not shown \label{fig:threeDTetSegmentA}] {
		\includegraphics[width=0.4 \textwidth]{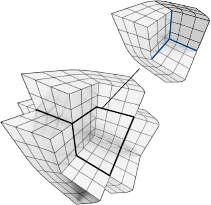} } 
	\hfill
	\subfloat[][Univariate weight functions in red and univariate B-splines in black and blue \label{fig:threeDTetSegmentB}] {
		\includegraphics[width=0.4\textwidth]{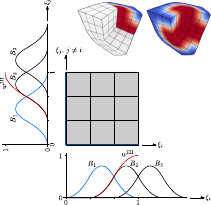} }
	\caption{Construction of the extraordinary joint weight function~${w^{J}_{\ell}}(\vec x)$ on one of the~$v = 4$ sectors. Each sector consists of the~$3 \times 3 \times 3$ elements shown in (a). The weight function~${w^{J}_{\ell}}(\vec \xi)$ is obtained using the univariate weight functions~$w^{\rm{III}}(\xi_1)$,~$w^{\rm{III}}(\xi_2)$ and~$w^{\rm{III}}(\xi_3)$ shown in (b) which are simply the complement to ones of~$w^{\rm{I}}(\xi_i) = 1 - w^{\rm{III}}(\xi_i)$ with~$i = 1, 2, 3$ shown earlier in Figure~\ref{fig:threeDPrismSegmentB}. The result is shown at the top right of (b) with the scalar field ranging between~$0$ (blue) and~$1$ (red). \label{fig:threeDTetSegment}}
\end{figure}

%
Without loss of generality, we consider a single extraordinary joint with valence~\mbox{$v = 4$} shared by four extraordinary prisms each of which have valence~\mbox{$e = 3$}, see Figure~\ref{fig:threeDSphere}, and describe the construction of its associated weight function. To simplify the construction of the extraordinary joint weight function~$w^J_{\ell} (\vec x)$ we require that the support of the already defined prism weight functions~${w^P_{k,\ell}}(\vec x)$ do not overlap at the joint. Recall from the Figures~\ref{fig:threeDPrismSegment} and~\ref{fig:threeDPrismWeight} that the support of the prism weight functions~${w^P_{k,\ell}}(\vec x)$ consist in the $\xi_1\xi_2$-plane corresponds to the $3$-neighbourhood of the extraordinary vertex. Hence, we choose the \mbox{$3$-neighbourhood} of the extraordinary vertex at the centre of the joint for constructing the weight function~${w^J_{\ell}}(\vec x)$ as illustrated in  Figure~\ref{fig:threeDTetSegmentA}. At the boundary of the 3-neighbourhood, the joint meets different prisms; the intersection of the joint with each prism is composed of $3\times 3  \times e$ quadrilateral faces, with~$e$ the valence of the prism's extraordinary edge. We require that across this $3\times 3  \times e$ of faces, the values and derivatives of the weight function~$w^J_\ell{(\vec x)}$ respectively match those of the
unique prism weight function that is non-zero on all these faces. For instance, in the already discussed example in Figure~\ref{fig:threeDPrismSegmentB} the weight functions~$w_{k,1}^P(\vec x)$ and  $w_{k,3}^P(\vec x)$ corresponding to~$w_{{k,1}}^{\rm{II}} (\xi_3)$ and~$w_{{k,3}}^{\rm{II}} (\xi_3)$ have to smoothly connect to the respective joint weight functions.

%

\begin{figure}[t!]
	\centering
	\subfloat[][${w^{J}_\ell}(\vec x)$ \label{fig:threeDTetWeightwQ1}] {
		\includegraphics[width=0.25\textwidth]{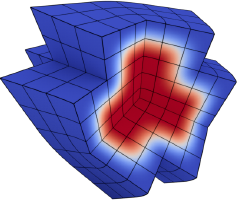} }
	\hfil
	\subfloat[][${w^P_{k,\ell}(\vec x)}$  \label{fig:threeDTetWeightwQk}] {
		\includegraphics[width=0.25 \textwidth]{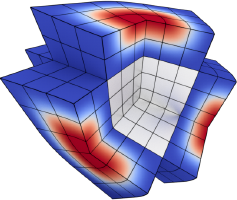} } 
	\hfil
	\subfloat[][$w^B(\vec x)$ \label{fig:threeDTetWeightwB}] {
		\includegraphics[width=0.25\textwidth]{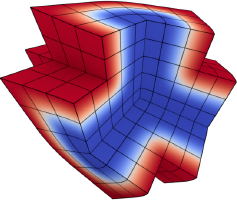} }
	\caption{Extraordinary joint weight function~${w^{J}_\ell}(\vec x)$, extraordinary prism weight functions~${w^P_{k,\ell}(\vec x)}$ and their complement to one $w^B(\vec x)$. The scalar field ranges between~$0$ (blue) and~$1$ (red). \label{fig:threeDTetWeightB}}
\end{figure}

To construct the joint weight function~${w^J_{\ell}}(\vec x)$ we follow once more a sector-wise approach as outlined in Figure~\ref{fig:threeDTetSegment}.  The joint is partitioned  into~$v = 4$ sectors each consisting of~$3 \times 3 \times 3$ elements. Each sector is parametrised using~\mbox{$\vec \xi = (\xi_1,\xi_2,\xi_3) \in \hat{\Omega}$} on the parametric domain \mbox{$\hat \Omega :=[0,1] \times [0,1] \times [0,1]$}. The extraordinary edges of the attached four prisms meet at the origin~$\vec \xi = (0,0,0)$ of  the domain~$\hat\Omega$. We define the joint weight function~$w_{\ell}^J(\vec \xi)$ as
\begin{align}
\begin{split}
	w_{\ell}^J(\vec \xi) = 1 & -  w^{\rm{III}}(\xi_1) \otimes w^{\rm{III}}(\xi_2) \otimes w^{\rm{III}}(\xi_3)  
	  -  \left( 1 -w^{\rm{III}}(\xi_1) \right )\otimes w^{\rm{III}}(\xi_2) \otimes w^{\rm{III}}(\xi_3) \\
	 & - w^{\rm{III}}(\xi_1) \otimes \left( 1 -w^{\rm{III}}(\xi_2) \right ) \otimes w^{\rm{III}}(\xi_3) 
	  -  w^{\rm{III}}(\xi_1) \otimes w^{\rm{III}}(\xi_2)  \otimes \left( 1 -w^{\rm{III}}(\xi_3) \right )  \,, 
\end{split}	 \label{eq:tetWeightFunction}
\end{align}
where the univariate weight functions~$w^{\rm{III}}(\xi_i)$ are, as depicted in Figure~\ref{fig:threeDTetSegmentB}, assembled from the smooth quadratic B-splines defined along the~$\xi_1$, $\xi_2$ or~$\xi_3$ axes. That is, 
\begin{equation} 
	w^{\rm{III}}(\xi_i) =  B_2(\xi_i) + B_3(\xi_i)  \quad i=1, \, 2, \, 3  \,.
\end{equation}
In Figure~\ref{fig:threeDTetWeightB} the isocontours of the  obtained joint weight function~${w^J_{\ell}}(\vec x)$,  the weight functions~${w^P_{k,\ell}(\vec x)}$ of the attached four prisms and their complement to one~$w^B(\vec x)$ are shown. In~\ref{sec:joint} we briefly demonstrate that the joint and prism weight functions for arbitrary vertex and edge valences~$v$ and~$e$ can be constructed following same approach.

%
\section{Examples \label{sec:examples}}
%
We proceed to establish the finite element convergence properties and accuracy of the SB-splines in solving Poisson and biharmonic problems. The respective weak forms and the details of the finite element discretisation are summarised in~\ref{sec:fem}. In all the examples we use sufficiently smooth manufactured solutions and focus on quadratic basis functions, except in 1D where we also consider cubic basis functions. As known from the isogeometric analysis literature the optimal converge rates for the Poisson problem discretised with standard quadratic B-splines are $3$ and $2$ in  the~$L^2$ and~$H^1$ (semi-)norms, respectively~\cite{Bazilevs:2006ab}. In contrast, the optimal convergence rates for the biharmonic problem discretised with quadratic  B-splines are $2$  in the  $L^2$ and the~$H^1$ (semi-)norms and $1$ in the ~$H^2$ (semi-)norm~\cite{tagliabue2014isogeometric}.
%
\subsection{One-dimensional Poisson problem} \label{sec:exampleOneD}
%
As a first example we consider the solution of a one-dimensional Poisson-Dirichlet problem~$- \D^2 u/ \D x^2 = f$ on the domain~$\Omega = (0,1)$. The body force~$f(x)$ is chosen such that the solution is
\begin{equation}
u(x) = \sin \left(3 \pi x \right) \,.
\end{equation}
The domain is parametrised using non-uniform B-splines of degree~$p_B=2$ and~$p_B=3$ in turn. The knot vector and control points~$x_i$ are selected so that each element has the same size~$h$. We intentionally introduce a~$C^0$ continuous kink at the midspan~$x = 1 / 2$ by setting the knot multiplicity to~$p_B$ therein. In addition, we use an open knot vector which allows the Dirichlet boundary condition to be imposed strongly.

\begin{figure}[h]
	\centering
	\subfloat[][Relative~$L^2$-norm of error\label{fig:1dPoissonConvergenceA}] {
		\includegraphics[width=0.32 \textwidth]{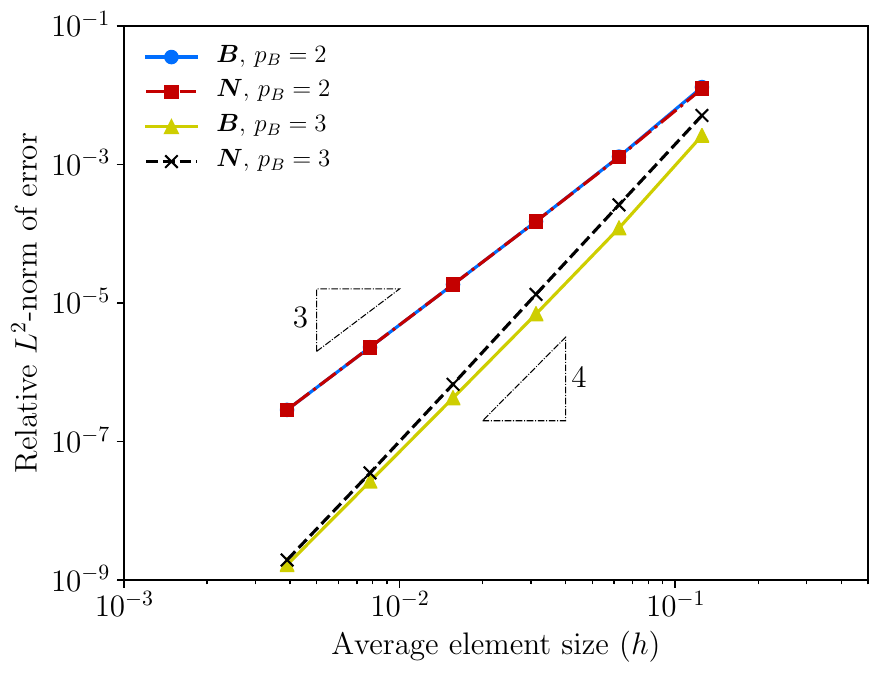} } 
	\hfill
	\subfloat[][Relative~$H^1$-seminorm of error\label{fig:1dPoissonConvergenceB}] {
		\includegraphics[width=0.32 \textwidth]{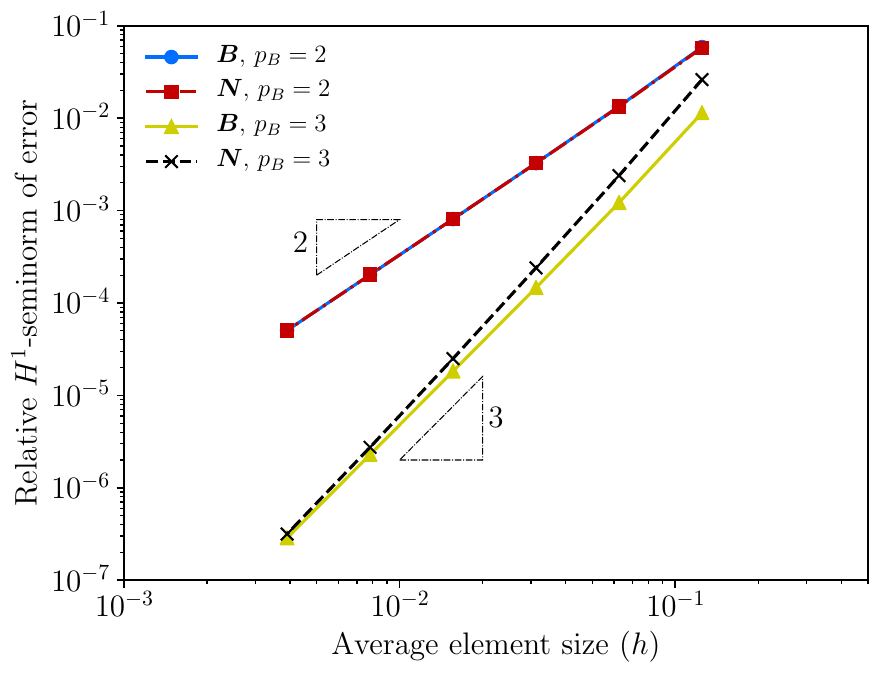} }
	\hfill
	\subfloat[][Condition number\label{fig:1dPoissonConvergenceC}] {
		\includegraphics[width=0.32 \textwidth]{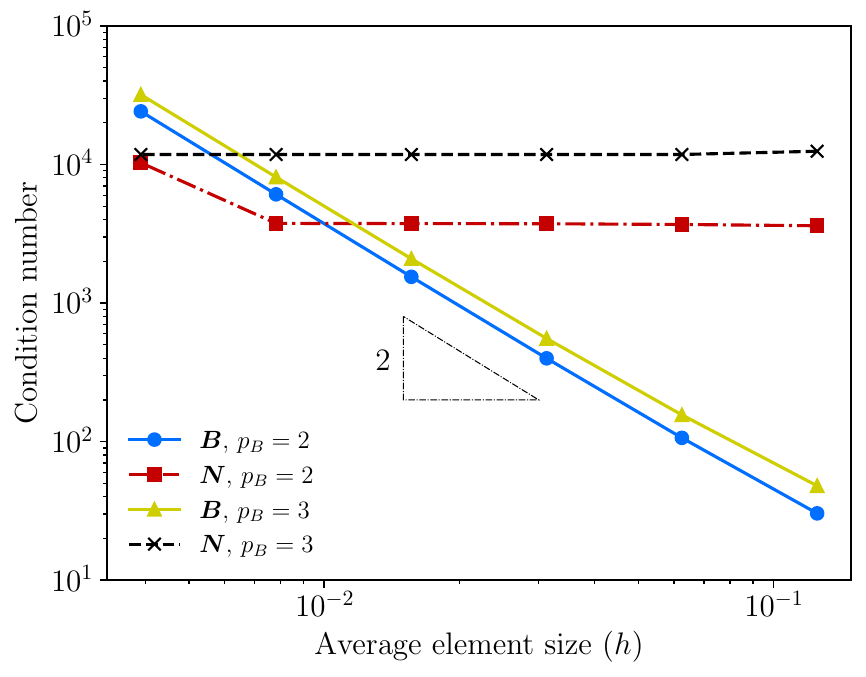} }
	\caption{One-dimensional Poisson-Dirichlet problem. Finite element convergence and condition number of the stiffness matrices for the $C^1$ continuous SB-splines~$\vec N(x)$ and the mixed smoothness B-splines~$\vec B(x)$, consisting of~$C^0$ and~$C^{p_B - 1}$ continuous basis functions.  \label{fig:1dPoissonConvergence}}
\end{figure}

In the following we compare the finite element convergence and the condition number of the stiffness matrices for the SB-splines~$\vec N(x)$ with the ones for the mixed B-splines~$\vec B(x)$, consisting of~$C^0$ and~$C^{p_B - 1}$ continuous basis functions. The $C^1$ continuous SB-splines~$\vec N(x)$ are constructed by blending B-splines~$\vec B(x)$ with Bernstein basis~$\vec Q(x)$ of same polynomial degree~$p_Q = p_B$. In comparison to B-splines~$\vec B(x)$ the SB-splines~$\vec N(x)$ have the additional degrees of freedom~$n_Q = p_Q + 1$. We begin with an initial coarse mesh of~$n_e = 8$ elements and obtain finer meshes using knot insertion. Figures~\ref{fig:1dPoissonConvergenceA} and~\ref{fig:1dPoissonConvergenceB} show that the SB-splines~$\vec N(x)$ yield optimal convergence rates for both polynomial degrees~$p_B = 2$ and~$p_B = 3$. In addition, the approximation error remains in the same order of magnitude with or without blending. The condition number of the respective stiffness matrices is plotted in Figure~\ref{fig:1dPoissonConvergenceC}. When the mesh size~$h$ is relatively large the SB-splines lead to significantly larger condition numbers than the B-splines. Interestingly, the condition numbers for SB-splines are almost  independent of mesh size. Overall, the SB-splines condition numbers compare favourably with the B-splines condition numbers.
%
\subsection{Poisson and biharmonic problems on a square domain} \label{sec:exampleSquare}
\begin{figure}[]
	\centering
		\includegraphics[width=0.4 \textwidth]{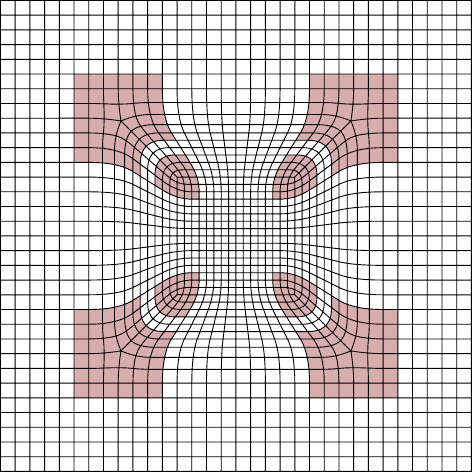} 
	\caption{Initial semi-structured coarse mesh of the square domain. The $32$ elements in the eight blending domains are shown shaded.  \label{fig:2dSquare} }
\end{figure}

%
We consider next the Poisson-Dirichlet and biharmonic problems on a square domain~\mbox{$\Omega = (0,1) \times (0, 1)$}. Figure~\ref{fig:2dSquare} shows the initial semi-structured coarse mesh with~$8$ extraordinary vertices. The boundary of the square domain is parametrised using open, uniform bi-quadratic B-splines. For the Poisson-Dirichlet problem, the Dirichlet boundary condition is imposed using Nitsche's method with the stabilisation parameter chosen as~$\gamma = 10/h^2$, see~\ref{sec:fem}. For the biharmonic problem, we use the penalty approach with the stabilisation parameter chosen as~$\gamma = 1/h^2$. In comparison to mixed B-splines~$\vec B(\vec x)$, for the SB-splines~$\vec N(\vec x)$ the additional degrees of freedom are~$n_Q = 8 \times 9 = 72$. We refine the mesh using a refinement scheme described in~\ref{sec:refinement}, so that the number of extraordinary vertices remains constant and the blending domains become increasingly smaller. In all meshes there are~$8$ extraordinary vertices and in total $32$ elements in the respective blending domains.

We approximate the finite element integrals using the Gauss-Legendre quadrature rule. In order to examine the effect of the number of quadrature points~$n_{gp}$ on the finite element convergence, for the SB-splines~$\vec N(\vec x)$, we vary~$n_{gp}$ for the domain integrals and use always~$3$ quadrature points for the boundary integrals. For the mixed B-splines~$\vec B(\vec x)$, we use for the domain integrals and boundary integrals~$3 \times 3$ and~$3$ quadrature points, respectively.

\begin{figure}[]
	\centering
	\subfloat[][Relative~$L^2$-norm of error] {
		\includegraphics[width=0.32 \textwidth]{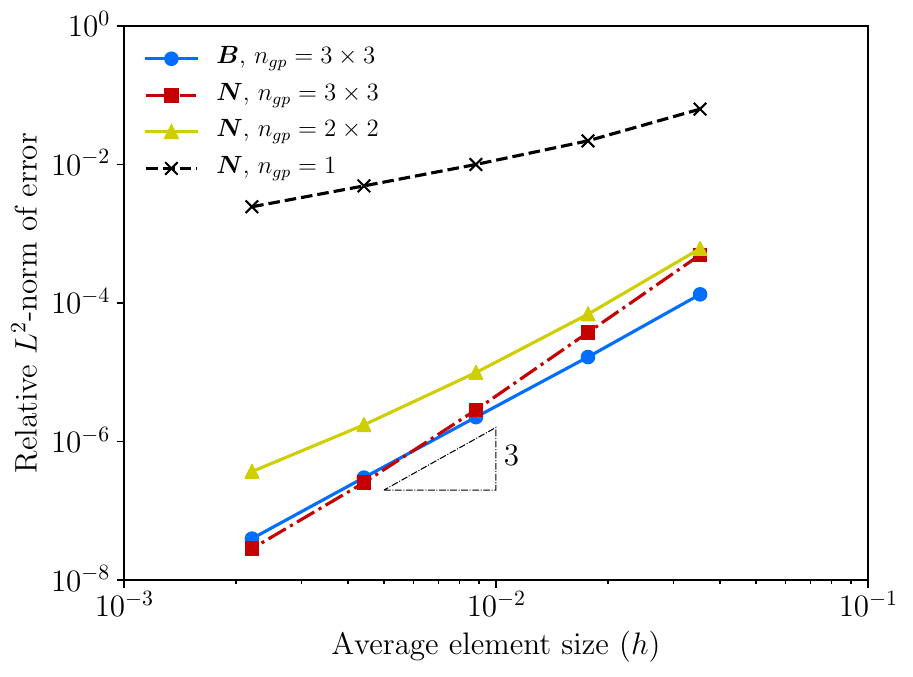} }
	\hfil
	\subfloat[][Relative~$H^1$-seminorm of error] {
		\includegraphics[width=0.32 \textwidth]{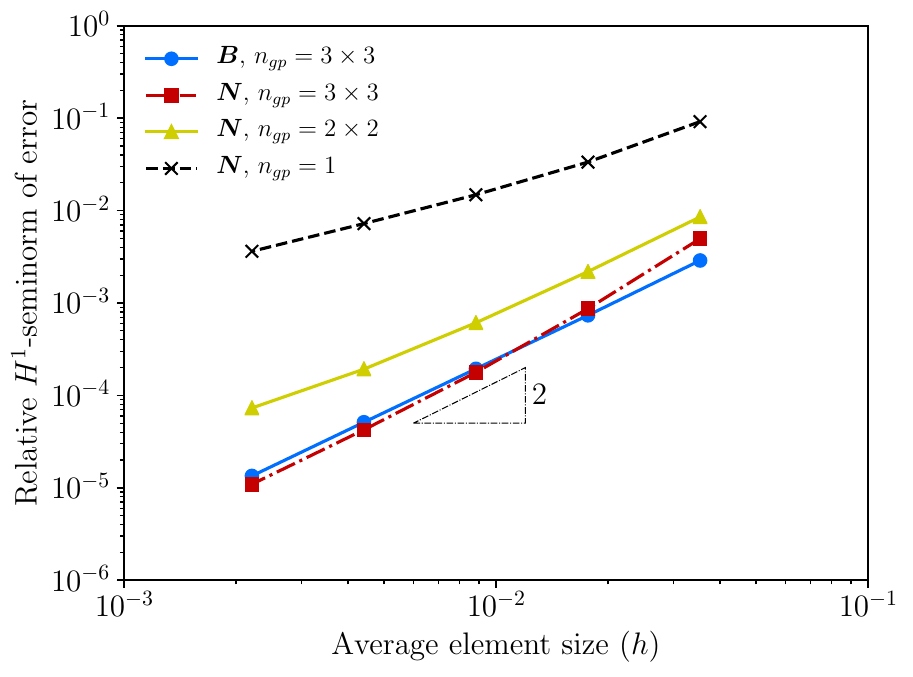} }
	\caption{Poisson-Dirichlet problem on a square domain. Convergence of the mixed B-spline~$\vec B(\vec x)$ and SB-spline~$\vec N(\vec x)$ solutions. \label{fig:2dPoissonSquare}}
\end{figure}

For the two-dimensional Poisson-Dirichlet problem, the body force~$f(\vec x)$ is chosen so that the solution is equal to
\begin{equation}
u(\vec x) = \sin\left( 6 x_1 \right)  \sin\left( 8 x_2 \right) \,.
\end{equation}
Figure~\ref{fig:2dPoissonSquare} confirms that the SB-splines are optimally convergent provided that at least~$n_{gp} = 3\times 3$ quadrature points are used. For~$n_{gp} = 3 \times 3$, as the mesh is refined the approximation error remains in the same order of magnitude with or without blending. We conjecture that the number of quadrature points~$n_{gp}$ for the SB-splines~$\vec N(\vec x)$ to achieve the optimal convergence rate is relatively small because the weight functions~$w^B(\vec x)$ and~$w^Q(\vec x)$ are assembled from smooth piecewise quadratic B-splines.

\begin{figure}[!t]
	\centering
	\subfloat[][Relative~$L^2$-norm of error] {
		\includegraphics[width=0.32 \textwidth]{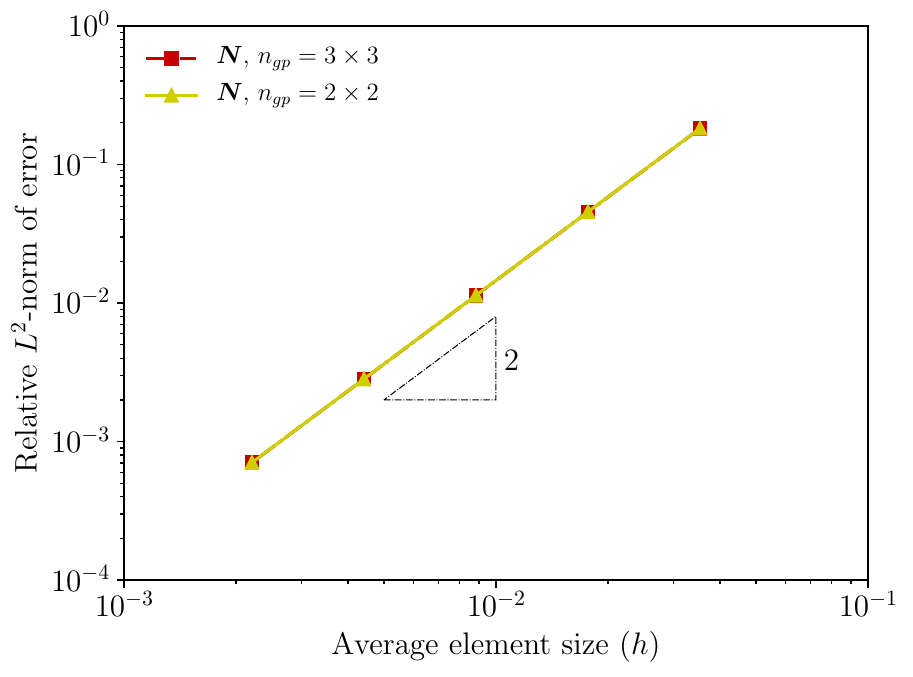} } 
	\hfil
	\subfloat[][Relative~$H^1$-seminorm of error] {
		\includegraphics[width=0.32 \textwidth]{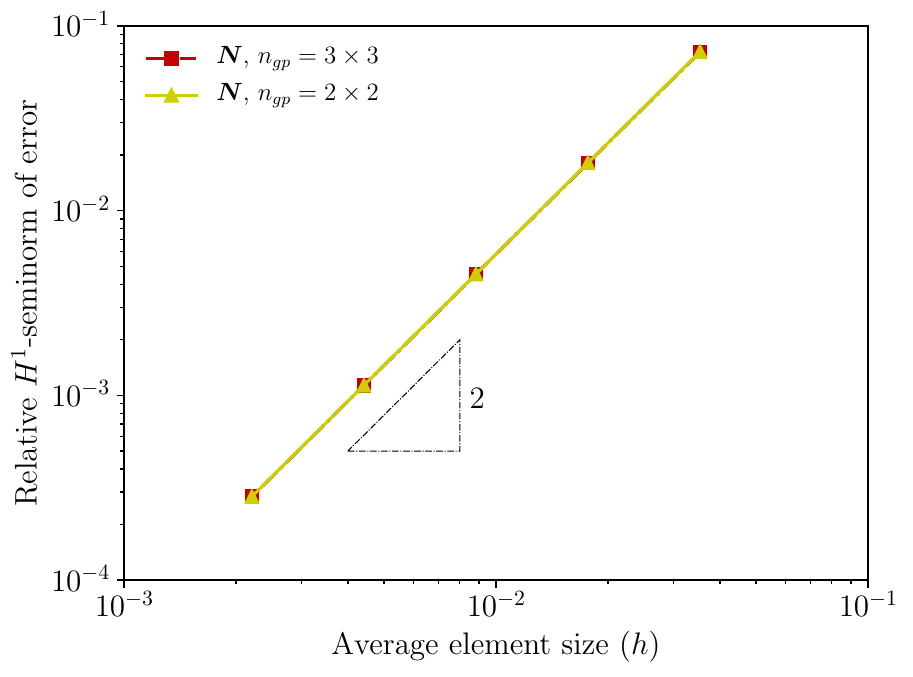} }
	\\ 
	\subfloat[][Relative~$H^2$-seminorm of error] {
		\includegraphics[width=0.32 \textwidth]{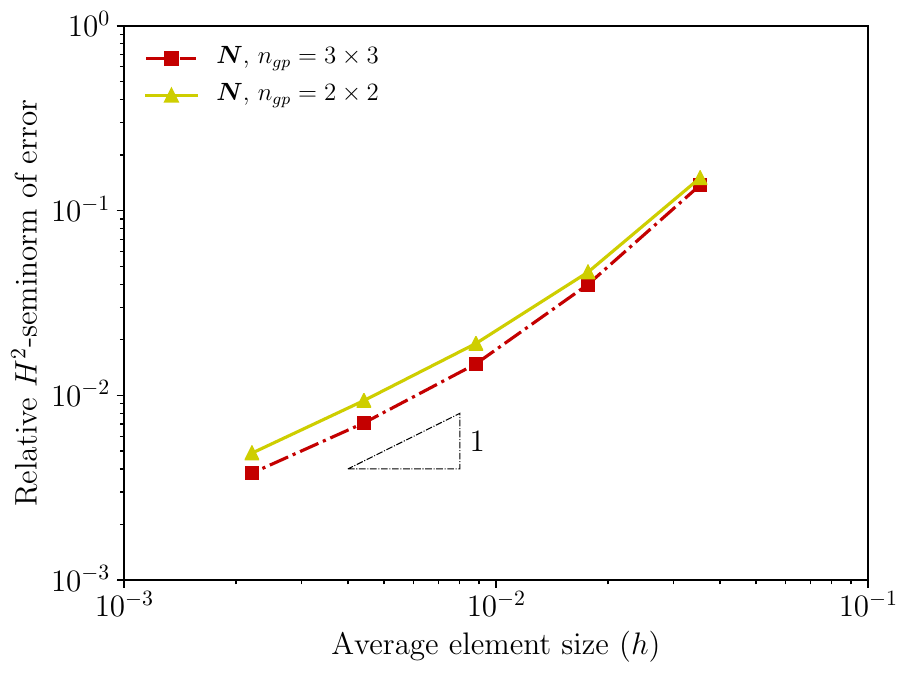} }
	\hfil
	\subfloat[][$L^\infty$-norm of error\label{fig:2dBiharmonicSquareMaxNormB}] {
		\includegraphics[width=0.32 \textwidth]{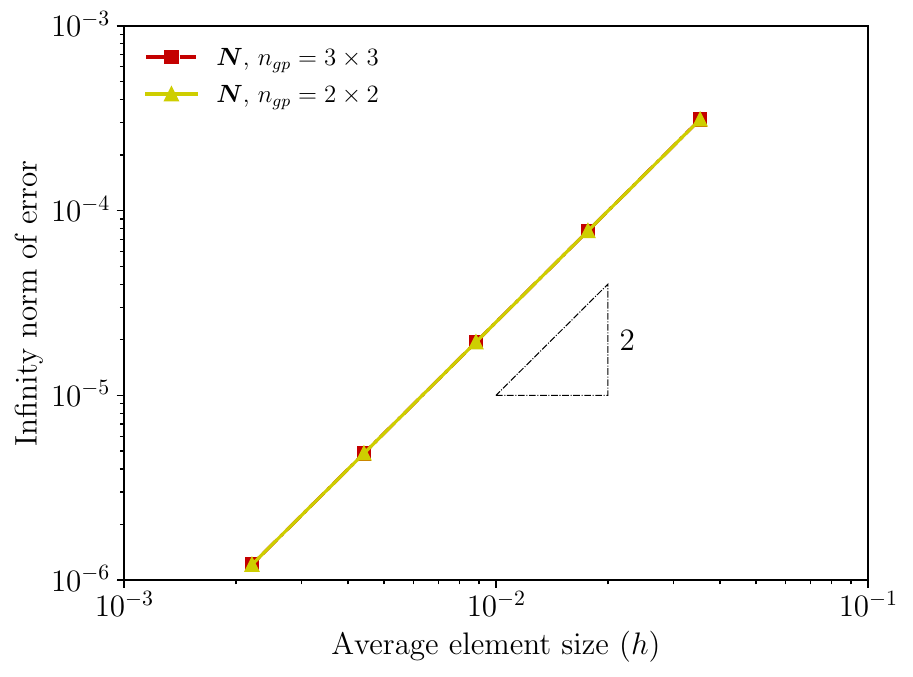} }	
	\caption{Biharmonic problem on a square domain. Convergence with SB-splines~$\vec N(\vec x)$. \label{fig:2dBiharmonicSquare}}
	
	\bigskip
	
	\subfloat[][$u^h$, out-of-plane scale factor: $100$\label{fig:2dBiharmonicField}] {
	\includegraphics[width=0.3 \textwidth]{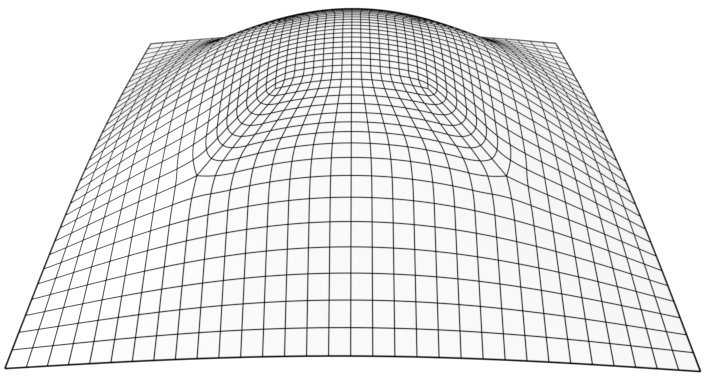} }
	\hfill
	\subfloat[][$\frac{\partial u^h}{\partial x_1}$, out-of-plane scale factor: $25$\label{fig:2dBiharmonicFirstDerivative}] {
		\includegraphics[width=0.3 \textwidth]{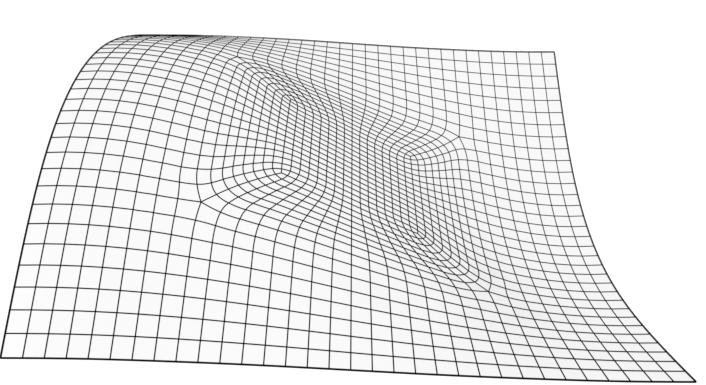} }
	\hfill
	\subfloat[][$\frac{\partial^2 u^h}{\partial x^2_1}$, out-of-plane scale factor: $2$\label{fig:2dBiharmonicSecondDerivative}] {
		\includegraphics[width=0.3 \textwidth]{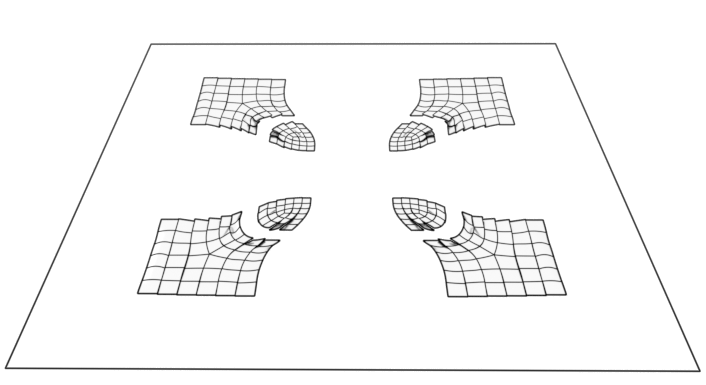} }
	\caption{Biharmonic problem on a square domain. The finite element solution~$u^h$ and its first and second partial derivatives with respect to~$x_1$ for the initial coarse mesh. All plots have been warped in the out-of-plane direction using the stated scale factors. In (a) and (b) all elements are shown whereas in (c) only elements in the blending domain~$\Omega_Q$ are shown. \label{fig:2dBiharmonicSquareWarp}}
\end{figure}

For the two-dimensional biharmonic problem, the body force~$f(\vec x)$ is chosen so that the solution is equal to
\begin{equation}
u(\vec x) = \frac{\sin\left( \pi x_1 \right)  \sin\left( \pi x_2 \right) }{4\pi^4}   \,.
\end{equation}
As the mixed B-splines~$\vec B(\vec x)$ are not globally~$C^1$ continuous on the considered semi-structured quadrilateral mesh, we examine only the finite element convergence using the~$C^1$ continuous SB-splines~$\vec N(\vec x)$. Figure~\ref{fig:2dBiharmonicSquare} shows that the SB-splines~$\vec N(\vec x)$ are optimally convergent for the biharmonic problem provided that a minimum of~$n_{gp} = 2 \times 2$ is used for the quadrature. However, note that the relative~$H^2$-seminorm of error is improved using~$n_{gp}= 3 \times 3$. 

In addition to the convergence rate, we examine the finite element solution~$u^h(\vec x)$ for the biharmonic problem qualitatively. Figure~\ref{fig:2dBiharmonicSquareWarp} shows the finite element solution~$u^h(\vec x)$ and its first and second partial derivatives with respect to~$x_1$ for the initial coarse mesh. Since the SB-splines~$\vec N(\vec x)$ are globally~$C^1$ continuous, both the finite element solution~$u^h(\vec x)$ and its first partial derivative with respect to~$x_1$ are continuous as visible in Figures~\ref{fig:2dBiharmonicField} and~\ref{fig:2dBiharmonicFirstDerivative}, respectively. Furthermore, as known the spatial derivatives of~$u^h(\vec x)$ often exhibit short-wavelength oscillations near the extraordinary vertices. Similarly, we observe such oscillations specifically for the second spatial derivatives in the blending domain~$\Omega_Q$ as shown in Figure~\ref{fig:2dBiharmonicSecondDerivative}. However, there is no oscillation in the~$1$-neighbourhood of the extraordinary vertices.
\begin{figure}[]
	\centering

\end{figure}
%
\subsection{Biharmonic problem on $v$-gon domains} \label{sec:exampleVGon}
\begin{figure}[]
	\centering
	\subfloat[][$v = 3$] {
		\includegraphics[width=0.32 \textwidth]{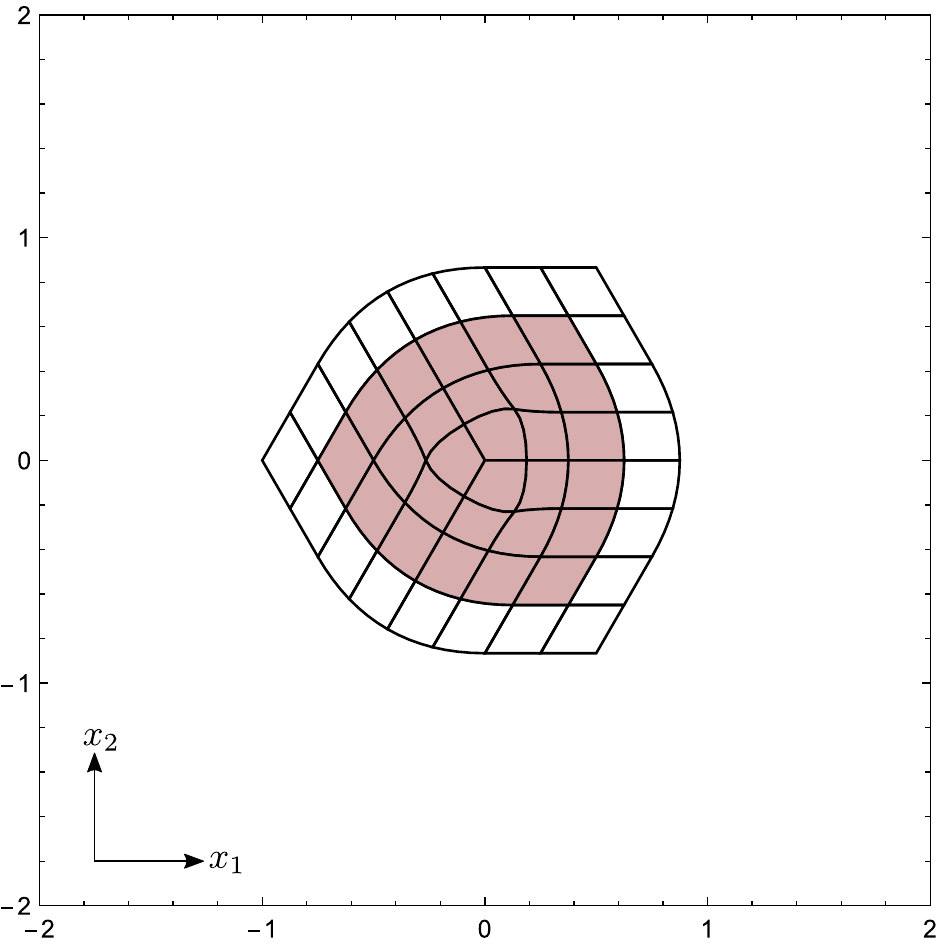} } 
	\hfill
	\subfloat[][$v = 5$] {
		\includegraphics[width=0.32 \textwidth]{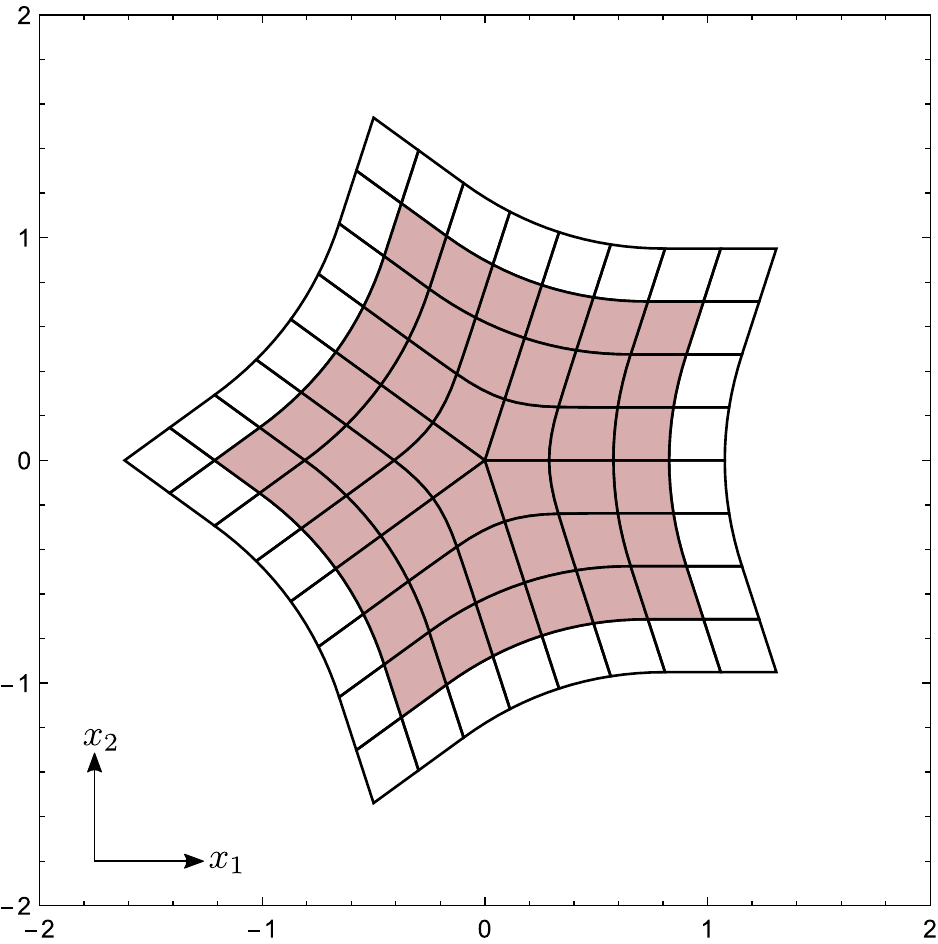} }
	\hfill
	\subfloat[][$v = 6$] {
		\includegraphics[width=0.32 \textwidth]{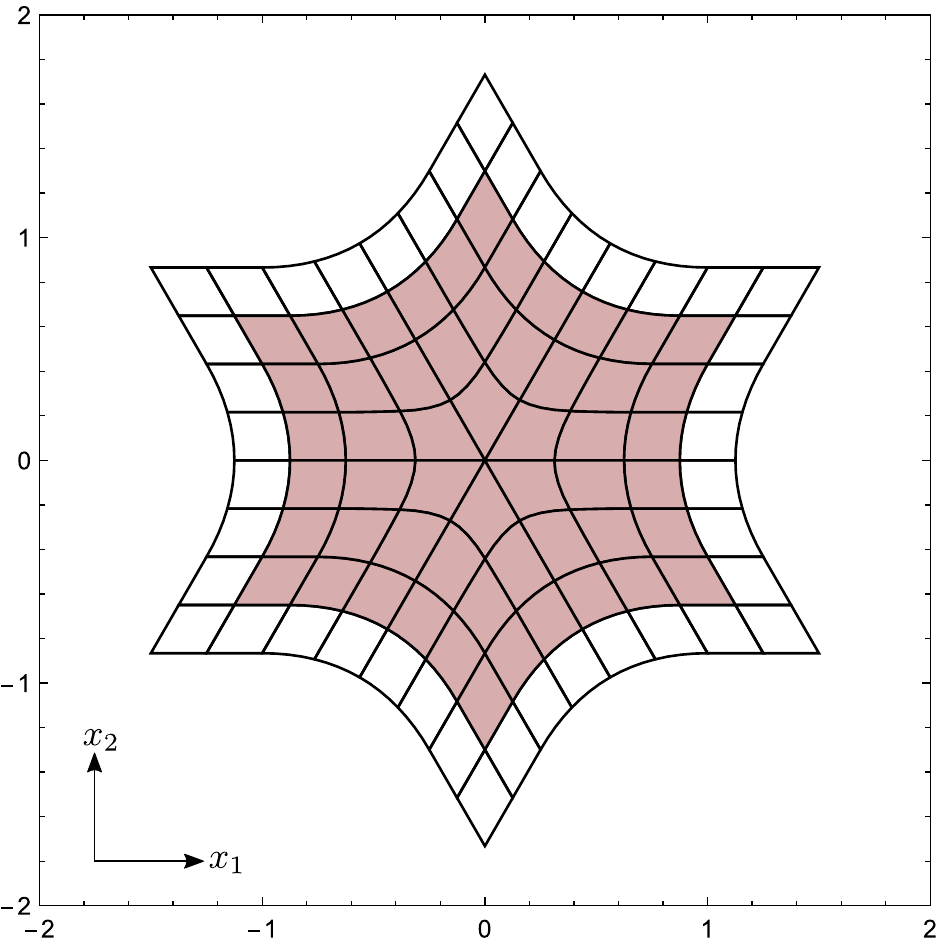} }
	\hfil
	\subfloat[][$v = 7$] {
		\includegraphics[width=0.32 \textwidth]{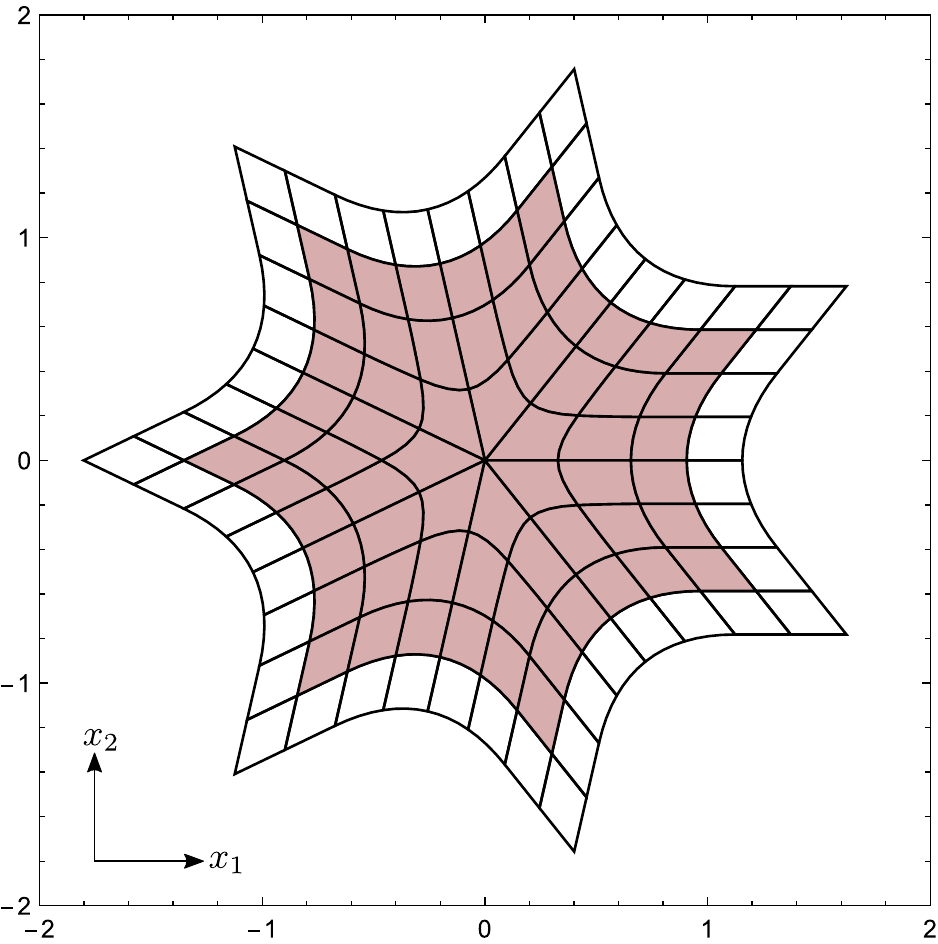} }
	\hfil
	\subfloat[][$v = 8$] {
		\includegraphics[width=0.32 \textwidth]{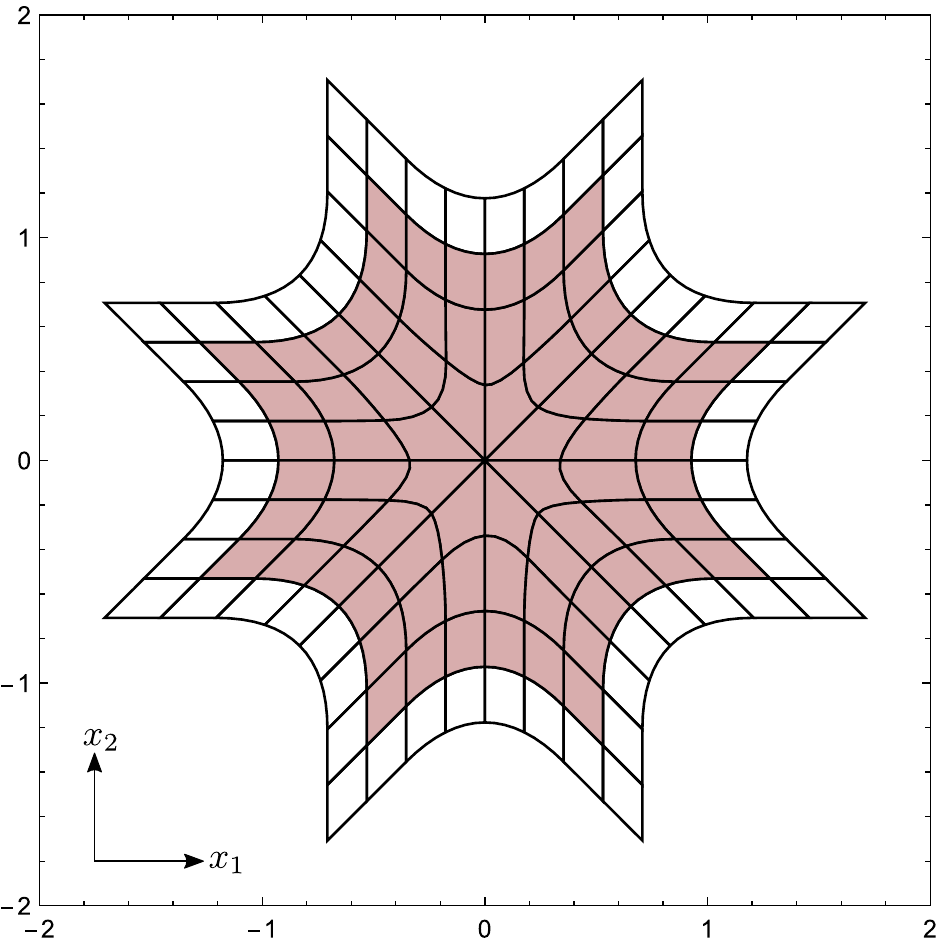} }
	\caption{Initial coarse meshes for the $v$-gon domains. For each initial coarse mesh the blending domain~$\Omega_Q$ is shaded in pink. \label{fig:2dvGon}}
\end{figure}

%
In some smooth basis function construction techniques, the respective finite element convergence rates are known to deteriorate when the valence~$v$ is increased, see the discussion in~\cite{toshniwal2017smooth}. Therefore, we investigate next the convergence rate of the SB-splines~$\vec N(\vec x)$ for different valences. To this end, we consider the biharmonic problem on five~$v$-gon domains with~\mbox{$v \in \{3, 5, 6, 7,8 \}$} as depicted in  Figure~\ref{fig:2dvGon}.  For each domain, the extraordinary vertex is located at the global origin~\mbox{$\vec x = (0 \quad 0)^\trans$}.
The boundary is parametrised with open, uniform bi-quadratic B-splines. To impose the boundary conditions, we use the penalty approach with the stabilisation parameter chosen as~$\gamma = 1000/h^2$. As shown to be sufficient for the biharmonic problem in Section~\ref{sec:exampleSquare}, we use $3 \times 3$ and~$3$ quadrature points for approximating the domain and boundary integrals, respectively. Similarly, we refine the mesh using the refinement scheme described in~\ref{sec:refinement}. The body force~$f(\vec x)$ is chosen such that the solution is
\begin{equation}
u( \vec x) = \sin \left(3 x_1 \right) \cos \left(3 x_2 \right) \,.
\end{equation}

\begin{figure}[]
	\centering
	\subfloat[][Relative~$L^2$-norm of error] {
		\includegraphics[width=0.32 \textwidth]{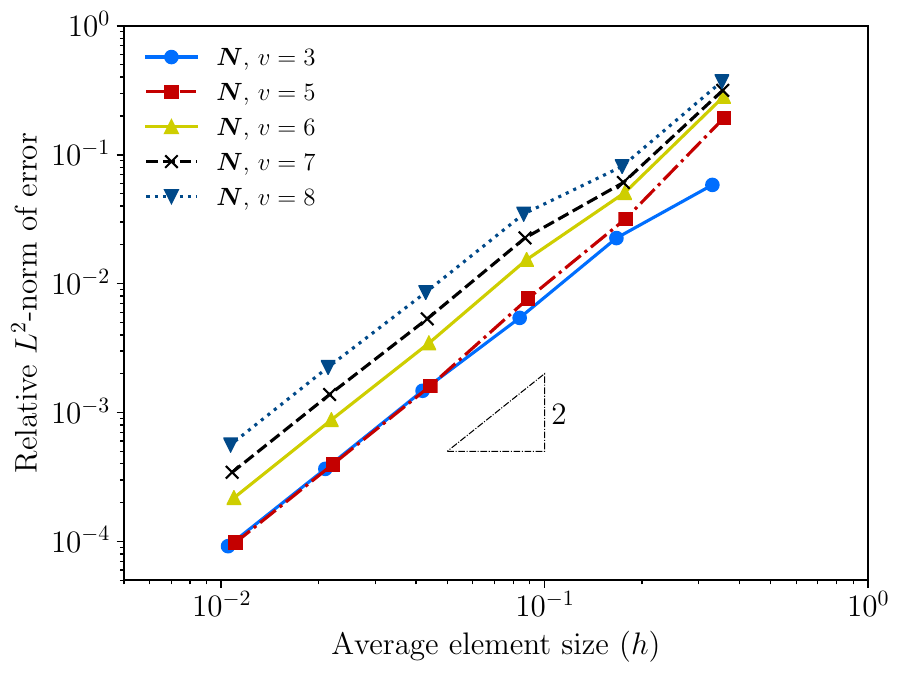} } 
	\hfill
	\subfloat[][Relative~$H^1$-seminorm of error] {
		\includegraphics[width=0.32 \textwidth]{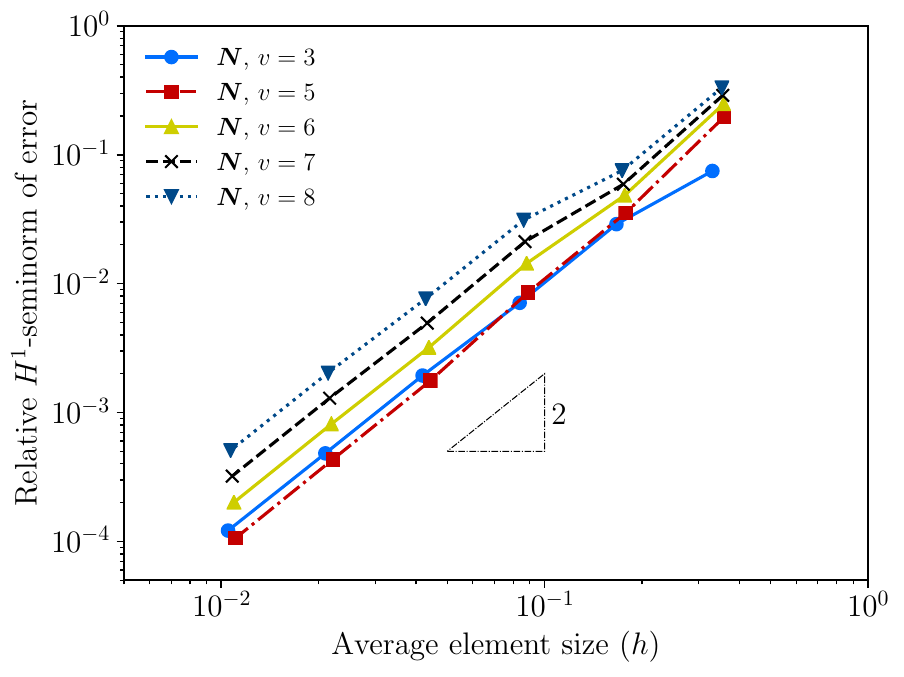} }
	\hfill
	\subfloat[][Relative~$H^2$-seminorm of error] {
		\includegraphics[width=0.32 \textwidth]{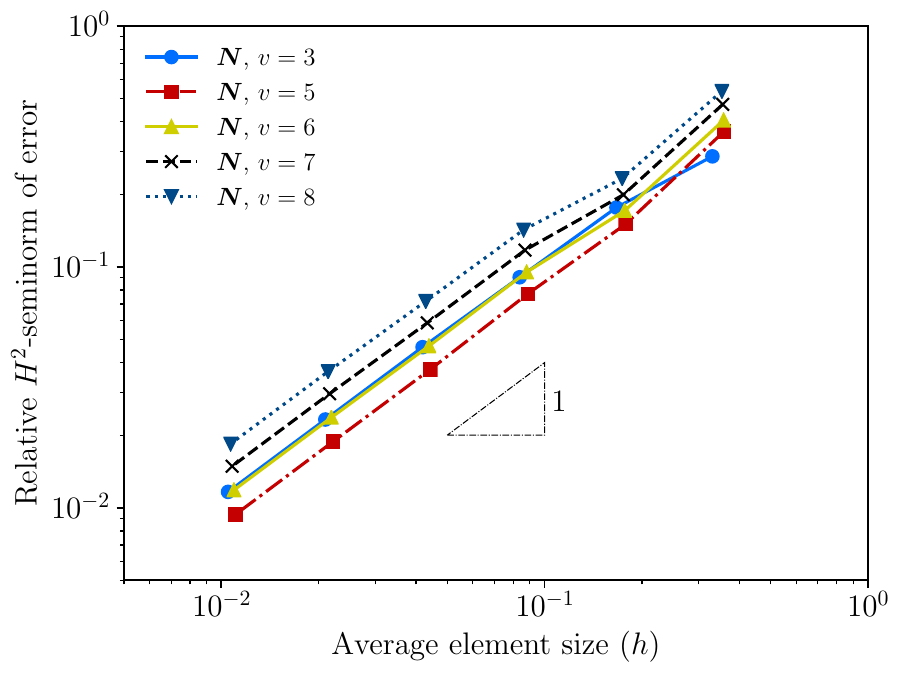} }
	\caption{Biharmonic problem on $v$-gon domains. Convergence with SB-splines~$\vec N(\vec x)$. \label{fig:2dBiharmonicVGon}}
\end{figure}

Figure~\ref{fig:2dBiharmonicVGon} shows the finite element convergence using the SB-splines~$\vec N(\vec x)$. Although the convergence rates for the first few coarser meshes are slightly fluctuating, overall the SB-splines~$\vec N(\vec x)$ are optimally convergent. In other words, the convergence rate is identical for the considered valences~$v \in \{3, 5, 6, 7,8 \}$. The increase of the valence leads, however, to a small increase in the convergence constants. This finding suggests that the SB-splines~$\vec N(\vec x)$ are robust since the studied valences~$v \in \{3, 5, 6, 7,8 \}$ are the most prevalent in well-designed meshes.
%
\subsection{Poisson and biharmonic problems on a spherical domain} \label{sec:exampleSphere}
%
As a final example, we consider the Poisson-Dirichlet and biharmonic problems on the spherical domain in Figure~\ref{fig:introSphere}. The spherical domain has a radius of~$2.55$ and is centred at the global origin~$\vec x = ( 0, 0, 0 )$. As an approximation to the spherical domain, the parametrised hexahedral mesh in Figure~\ref{fig:introSphereParam} has an average mesh size~$h = 0.4228$ and consists of~$20$ extraordinary prisms of valence~$e = 3$ and~$8$ extraordinary joints. The boundary of the spherical domain is parametrised using open, uniform tri-quadratic B-splines. For the Poisson-Dirichlet problem, the Dirichlet boundary condition is imposed using Nitsche's method with the stabilisation parameter chosen as~$\gamma = 10/h^2$. For the biharmonic problem, we use the penalty approach with the stabilisation parameter chosen as~$\gamma = 1000/h^2$. The numbers of degrees of freedom for the mixed B-splines~$\vec B(\vec x)$ and the SB-splines~$\vec N(\vec x)$ are~$6413$ and~$6413 + (20 + 8) \times 27 = 7169$, respectively.

\begin{figure}[]
	\centering
	\subfloat[][Analytical $\frac{\partial{u(\vec x)}}{\partial x_1} $ \label{fig:exampleSphereAnalytical}] {
		\includegraphics[width=0.3 \textwidth]{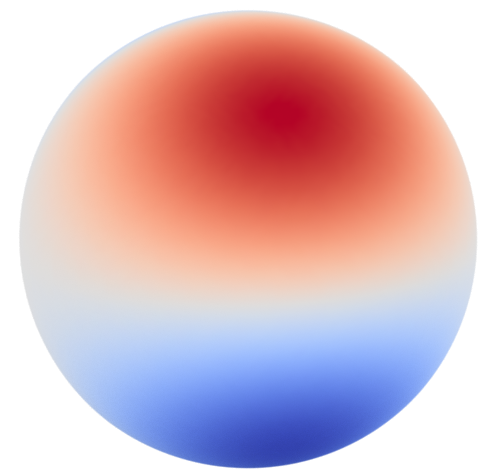} }
	\hfill
	\subfloat[][$\frac{\partial{u^h(\vec x)}}{\partial x_1}$ using~$\vec B(\vec x)$  \label{fig:exampleSphereBSplines}] {
		\includegraphics[width=0.3 \textwidth]{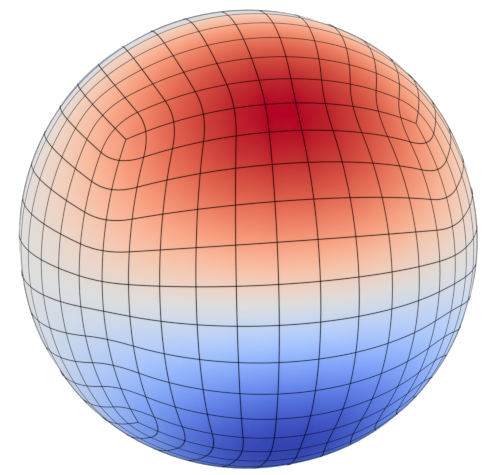} }
	\hfill
	\subfloat[][$\frac{\partial{u^h(\vec x)}}{\partial x_1}$ using~$\vec N(\vec x)$  \label{fig:exampleSphereBlended}] {
		\includegraphics[width=0.3 \textwidth]{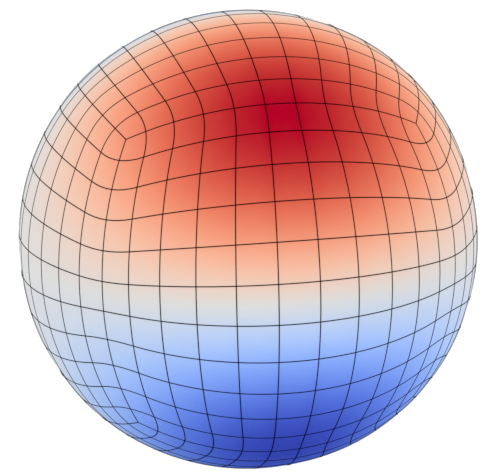} }
	\caption{Poisson-Dirichlet problem on a spherical domain. Isocontours of the first partial derivative of the field solution with respect to~$x_1$ on the spline surface. The scalar field ranges between $-0.18$ (blue) and $0.18$ (red).\label{fig:3dSphereDerivative}}
\end{figure}

\begin{figure}[]
	\centering
	\subfloat[][$\frac{\partial{u^h(\vec x)}}{\partial x_1}$ using~$\vec B(\vec x)$ \label{fig:exampleSphereC0}] {
		\includegraphics[width=0.4 \textwidth]{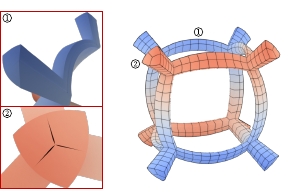} }
	\hfill
	\subfloat[][$\frac{\partial{u^h(\vec x)}}{\partial x_1}$ using~$\vec N(\vec x)$  \label{fig:exampleSphereC1}] {
		\includegraphics[width=0.4 \textwidth]{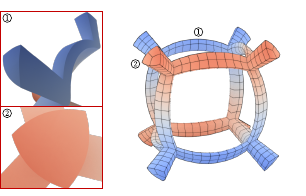} }
	\hfill
	\caption{Poisson-Dirichlet problem on a spherical domain. Isocontours of the first partial derivative of the finite element solution with respect to~$x_1$ near the extraordinary edges. For consistency, the scalar field herein has the same scale as Figure~\ref{fig:3dSphereDerivative}.\label{fig:3dSphereContinuity}}
\end{figure}

For the three-dimensional Poisson-Dirichlet problem, the body force~$f(\vec x)$ is chosen so that the solution is
\begin{equation}
u(\vec x) = \sin\left( \frac{x_1}{2}  \right)  \sin\left( \frac{x_2}{2} \right)  \sin\left( \frac{x_3}{4}  \right)   \,.
\end{equation}
We compare numerically the finite element solution~$u^h(\vec x)$ between~$\vec B(\vec x)$ and~$\vec N(\vec x)$. The relative $L^2$-norms of error for~$\vec B(\vec x)$ and~$\vec N(\vec x)$ are~$4.9111 \times 10^{-4}$ and $3.9034 \times 10^{-4}$, respectively whereas the relative~$H^1$-seminorms of error for~$\vec B(\vec x)$ and~$\vec N(\vec x)$ are $5.1337 \times \times 10^{-3}$ and~$3.5630 \times 10^{-3}$, respectively. Therefore, the approximation error has the same order of magnitude with and without blending. Figure~\ref{fig:3dSphereDerivative} illustrates the first partial derivative of the finite element solution with respect to~$x_1$. As shown in the figure, both~$\vec B(\vec x)$ and~$\vec N(\vec x)$ yield an accurate approximation to the first partial derivative of the analytical solution with respect to~$x_1$. In addition, Figure~\ref{fig:3dSphereContinuity} ascertains that the SB-splines~$\vec N(\vec x)$ are globally~$C^1$ continuous. For instance, the mixed B-splines~$\vec B(\vec x)$ yield~$C^0$ continuous finite element solution~$u^h(\vec x)$ near the extraordinary edges as inferred from the discontinuity of the first partial derivative of~$u^h$ with respect to~$x_1$ in Figure~\ref{fig:exampleSphereC0}. In contrast, the first partial derivative of~$u^h$ with respect to~$x_1$ using~$\vec N(\vec x)$ is continuous in the same subdomains as shown in Figure~\ref{fig:exampleSphereC1}.

\begin{figure}[t]
	\centering
	\subfloat[][Analytical $u(\vec x)$ \label{fig:exampleSphereBiharmonicAnalytical}] {
		\includegraphics[width=0.3 \textwidth]{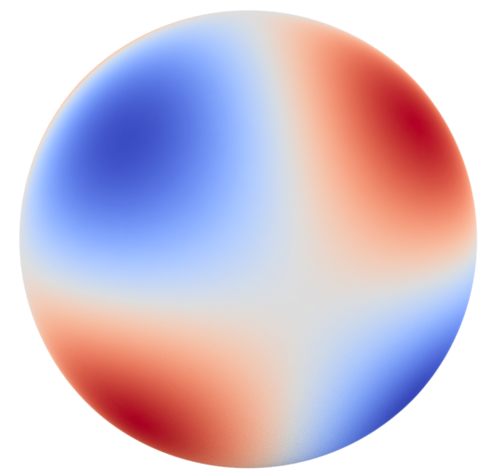} }
	\hfil
	\subfloat[][$u^h(\vec x)$ using~$\vec N(\vec x)$  \label{fig:exampleSphereBiharmonicBlended}] {
		\includegraphics[width=0.3 \textwidth]{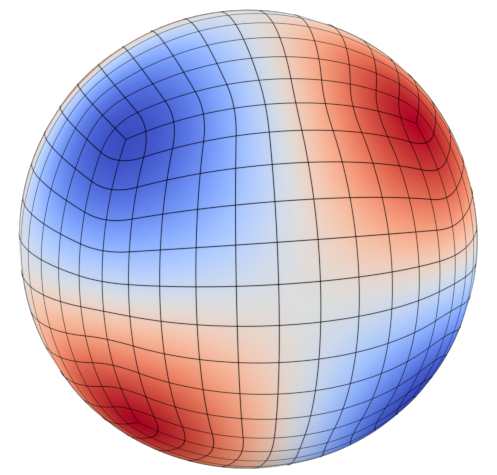} }
	\caption{Biharmonic problem on a spherical domain. Isocontours of the field solution on the spline surface. The scalar field ranges between $-0.43$ (blue) and $0.43$ (red).\label{fig:3dSphereBiharmonic}}
\end{figure}

For the three-dimensional biharmonic problem, the body force~$f(\vec x)$ is chosen so that the solution is equal to
\begin{equation}
u(\vec x) = \frac{1}{8 \pi} \sin\left( \frac{ \pi x_1}{8}  \right)  \sin\left( \frac{ \pi x_2}{8} \right)  \sin\left( \frac{\pi x_3}{8}  \right)    \,.
\end{equation}
Using the SB-splines~$\vec N(\vec x)$, the relative $L^2$-norm of error, relative~$H^1$-seminorm of error and relative~$H^2$-seminorm of error are~$9.8943 \times 10^{-4}$,~$3.7825 \times 10^{-3}$ and~$6.4669 \times 10^{-2}$, respectively.  Figure~\ref{fig:3dSphereBiharmonic} shows that the SB-splines~$\vec N(\vec x)$ yield a satisfactory finite element approximation to the analytical solution.
\section{Conclusions \label{sec:conclusions}}
%
We introduced SB-splines i.e.\ a smooth blended B-spline construction for unstructured quadrilateral and hexahedral meshes and demonstrated its optimal convergence for quadratic B-splines on unstructured quadrilateral and hexahedral meshes. We determine the smooth weight functions required for blending from the smooth mixed B-splines defined on the regular parts of the unstructured mesh. The weight functions multiplied with the available mixed B-splines and additionally introduced Bernstein basis functions yield new basis functions. As shown numerically, the SB-splines can be efficiently integrated using standard Gauss-Legendre quadrature with a very small number of quadrature points. In the blending region, the new basis functions have slightly larger support close to the extraordinary features. For instance, in 2D the support consists of the 3-neighbourhood of the extraordinary vertex. Remarkably, the numerically determined convergence rates in 2D are optimal for both Poisson and biharmonic problems and are independent of the valence of the extraordinary vertex. However, the convergence constants show a slight increase with an increase in valence, which may be explained by the short-wavelength oscillations, or ripples, in the higher-order derivatives at the blending region. As discussed, on unstructured hexahedral meshes, the extraordinary edges and vertices usually form a connected network. The respective weight functions obtained from the available smooth mixed B-splines may not have a compact support.  Therefore, we decompose the weight functions so that the resulting weight functions and basis functions have a compact support and are still polynomials in the parameter space.

In closing,  we stress that the proposed construction can be applied to mixed B-splines of arbitrary degrees, although we have studied only quadratic mixed B-splines so far. To this end, it is necessary to extend the introduced mesh refinement scheme for quadratic mixed B-splines to arbitrary degree. In our experience, the details of this refinement are important for achieving optimal convergence rates. For cubic SB-splines the mixed $C^0/C^1/C^2$ B-splines introduced in Wei et al.~\cite{wei2018blended} appear as particularly promising. Moreover, while we presented some mathematical analyses (e.g., proof of linear independence in 1D), further analysis is needed to prove the numerically observed properties of SB-splines. In principle, the proposed construction can also be applied to non-uniform B-splines and extended to surfaces with arbitrary topology, i.e. 2-manifolds in $\mathbb{R}^3$. Non-uniform constructions can, amongst others, significantly simplify the enforcement of essential boundary conditions. In case of surfaces with arbitrary topology, the SB-splines have to be constructed on a set of intermediate parametric domains corresponding to each of the extraordinary vertices.  The so obtained SB-splines on parametric domains can be subsequently mapped to~$\mathbb R^3$. Lastly, to make the presented construction useful in geometric design, the introduced new degrees of freedom around the extraordinary features must be associated with control vertex positions. Following related constructions in geometric modelling, see~\cite{ying2004simple,zorin2006constructing,levin2006modified,antonelli2013subdivision,majeedCirak:2016}, this may be achieved by projecting the new degrees of freedom to the existing or possibly some new control vertex positions in the mesh.

\appendix
\section{Linear independence \label{sec:independence}}
%
We provide in this appendix a proof for the linear independence of the one-dimensional SB-splines. We consider as in Figures~\ref{fig:oneDQuadratic} and~\ref{fig:oneDCubic} a 1D setup with $n_B$ mixed B-splines $B_i(x)$ of degree $p_B \geq 2$. As discussed in Section \ref{sec:oneD}, we assume that the prescribed B-spline smoothness is $C^0$ at a single breakpoint and $C^{p_B-1}$ at all others. In general, the mixed B-splines $B_i(x)$ will be non-polynomials in physical space because of the isoparametric mapping, see~\eqref{eq:oneDMappingBSplines}. Let $\mathcal{B}$ denote the set of~$n_B$ mixed B-splines $B_i(x)$, $\mathcal{P}$ the set of~$n_Q = p_Q + 1$ Bernstein polynomials and $\mathcal{N}$ the set of $n_N = n_B + n_Q$ SB-splines. The set $\mathcal{B}$ is split into the two non-intersecting sets
\begin{equation}
\mathcal{B}^B = \left\{ B_i \mid \supp B_i \subset \Omega^Q \right\} \quad \text{and}\quad \mathcal{B}^C = \mathcal{B} \setminus \mathcal{B}^B \,.
\end{equation}
The set $\mathcal{N}$ is composed of the three non-intersecting sets
\begin{equation}
\mathcal{N}^B = \left\{ w^B B_i \mid B_i \in \mathcal{B}^B \right\} \text{ ,} \quad \mathcal{N}^Q = \left\{ w^Q Q_i \mid Q_i \in \mathcal{P} \right \} \quad \text{and} \quad \mathcal{N}^C = \mathcal{N} \setminus \left( \mathcal{N}^B \cup \mathcal{N}^Q \right) \,.
\end{equation}
As implied by the choice of weight function in \eqref{eq:oneDWeightChoice}, $\nexists B_i \in \mathcal{B}^B$ such that $w^B|_{\supp B_i} = 0$, and similarly $\nexists Q_i \in \mathcal{P}$ such that $w^Q|_{\supp Q_i} = 0$. We want to prove that the functions in $\mathcal{N}$ are linearly independent.

Linear independence requires that
\begin{equation}
\sum_{N_i \, \in \, \mathcal{N}}^{} N_i(x) \alpha_i = 0 \quad \forall x \in \Omega \,,
\label{eq:proof1}
\end{equation}
is satisfied only when the coefficients~$\alpha_i = 0$. Observe that outside the blending domain,~\eqref{eq:proof1} reduces to
\begin{equation}
\sum_{N_i \, \in \, \mathcal{N}^C}^{} N_i(x) \alpha_i = \sum_{B_i \, \in \, \mathcal{B}^C}^{} B_i(x) \alpha_i = 0 \quad \forall x \in \Omega \setminus \Omega^Q \,.
\end{equation}
Therefore, due to the linear independence of B-splines, we obtain~$\alpha_i = 0$ for all $N_i(x) \in \mathcal{N}^C$. The remaining terms in~\eqref{eq:proof1} correspond to splines with a support inside the blending region~$\Omega^Q $. We prove by contradiction that the coefficients of the non-vanishing terms must be zeros. Assume that the SB-splines are linearly dependent such that
\begin{equation}
\sum_{N_i \, \in \, \mathcal{N}^B}^{} N_i(x) \alpha_i = \sum_{N_j \, \in \, \mathcal{N}^Q}^{} N_j(x) \alpha_j  \quad \forall x \in \Omega^Q \,,
\end{equation}
or equivalently that
\begin{equation}
w^B(x)   \sum_{B_i \, \in \, \mathcal{B}^B}^{}  B_i(x) \alpha_i  = \left(1 - w^B(x) \right) \sum_{Q_j \, \in \, \mathcal{P}}^{}  Q_j(x) \alpha_j  \,.
\label{eq:proof2}
\end{equation}
Next, observe that there are two elements $\Omega_{k1}, \Omega_{k2} \subset \Omega^Q$ such that on each there is only one function from $\mathcal{N}^B$ that is non-zero, e.g., the leftmost or the rightmost element in the grey region in Figure \ref{fig:oneDQuadratic} or Figure \ref{fig:oneDCubic}.
Let the corresponding non-zero functions be $N_{i1}, N_{i2}$, respectively.
Then, we have $w^B|_{\Omega_{k\ell}} = (1-B_{i\ell})$ and $w^Q|_{\Omega_{k\ell}} = B_{i\ell}$, $\ell = 1, 2$, thus
\begin{equation}
(1 - B_{i\ell}(x)) \alpha_{i\ell} = \sum_{Q_j \, \in \, \mathcal{P}}^{}  Q_j(x) \alpha_j \quad \forall x \in \Omega_{k\ell} \subset \Omega^Q \,,\quad \ell = 1,2\;.
\label{eq:proof3a}
\end{equation}
Note that the right hand side is a polynomial function.
\begin{itemize}
\item[--] Case 1: Let $B_{i\ell}|_{\Omega_{k\ell}}$ be non-polynomial.
Then, the equality in \eqref{eq:proof3a} can be satisfied only if both sides are equal to $0$.
\item[--] Case 2: Let the isoparametric mapping be such that both $\Omega_{k1}, \Omega_{k2}$ are obtained by affinely mapping the associated element in parameter space.
Thus, $B_{i\ell}|_{\Omega_{k\ell}}$ are degree $p_B$ polynomials for both $\ell = 1, 2$.
This has two implications.
First, for equality, both the left and right hand sides in \eqref{eq:proof3a} need to represent the same polynomial, say $f$, of degree $p = \min\{p_B,p_Q \}$.
Note that $f$ is thus a global polynomial on $\Omega^Q$.
Second, by the end-point vanishing property of B-splines, $B_{i\ell}$ vanishes $p_B$ times on one of the endpoints of $\Omega_{k\ell}$, $\ell = 1, 2$.
This imposes $2p_B$ constraints on the polynomial $f$, thus implying $f = 0$.
\end{itemize}
Both the above cases imply that the right hand side in \eqref{eq:proof3a} is zero and, in particular, all coefficients of $Q_i$ are thus zero by their linear independence. As a result, the right hand side in \eqref{eq:proof2} is zero, thus implying that all coefficients of $B_i \in \mathcal{B}^B$ are zero by their linear independence.
%
\section{Finite element discretisation \label{sec:fem}}
%
\subsection{Poisson equation}
%
The Poisson equation is given by
\begin{equation} \label{eq:femPoissonStrongForm}
\begin{split}
\begin{alignedat}{2}
- \Delta u & = f                 \qquad &&\text{in } \Omega \,,\\
u &= \bar{u}                     \qquad &&\text{on } \Gamma_D \,, \\
\nabla u \cdot \vec n   &= \bar{t} \qquad &&\text{on } \Gamma_N \,,
\end{alignedat}
\end{split}
\end{equation}
where~$u$ is the solution field in the domain~$\Omega$ due to the body force~$f$, $\bar{u}$ is the prescribed solution field on the Dirichlet boundary~$\Gamma_D$,~$\bar{t}$ is the prescribed flux on the Neumann boundary~$\Gamma_N$ with the outward unit normal~$\vec n$,~$\nabla$ is the gradient operator and~$\Delta = \nabla \cdot \nabla$ is the Laplacian operator. The weak formulation of the Poisson equation can be stated as~\cite{Nitsche:1971aa,Fernandez-Mendez:2004aa}: Find~$u \in H^1(\Omega)$ such that
\begin{equation} \label{eq:femPoissonWeakForm}
a( u, v ) =  l(v) \, ,
\end{equation}
for all~$v \in H^1(\Omega)$ with 
\begin{subequations}
	\begin{align}
	a(u,v) &= \int_\Omega \nabla{u} \cdot \nabla{v}\D\Omega + \gamma \int_{\Gamma_D} uv\D\Gamma -  \int_{\Gamma_D} \Big( u \left( \nabla v \cdot \vec n \right) + v \left(\nabla u \cdot \vec n \right) \Big) \D\Gamma \,, 
	\\ 
	l(v) & = \int_\Omega v s \D \Omega + \int_{\Gamma_N} v \bar{t}  \D\Gamma + \gamma \int_{\Gamma_D} v \bar{u} \D \Gamma - \int_{\Gamma_D} \left( \nabla v \cdot \vec n \right) \bar{u} \D\Gamma \,, 
	\end{align}
\end{subequations}
and the positive stabilisation parameter~$\gamma$.
%
\subsection{Biharmonic equation}
%
%
The biharmonic equation is given by
\begin{equation} \label{eq:femBiharmonicStrongForm}
\begin{split}
\begin{alignedat}{2}
\Delta^2 u & = f                 \qquad &&\text{in } \Omega \,,\\
u = \bar{u}, \quad \nabla u \cdot \vec n &= \bar{t}               \qquad &&\text{on } \Gamma_D \,, \\
\Delta u = \bar{\kappa}, \quad  \nabla (\Delta u) \cdot \vec n  &= \bar{\lambda} \qquad &&\text{on } \Gamma_N \,,
\end{alignedat}
\end{split}
\end{equation}
where~$\bar{\kappa}$ and~$\bar{\lambda}$ are respectively the bending moment and shear force prescribed on the Neumann boundary~$\Gamma_N$. The weak formulation of the biharmonic equation can be stated as~\cite{embar2010imposing}: Find~$u \in H^2(\Omega)$ such that
\begin{equation} \label{eq:femBiharmonicWeakForm}
a( u, v ) =  l(v) \, ,
\end{equation}
for all~$v \in H^2(\Omega)$ where
\begin{subequations}
	\begin{align}
	\nonumber 
	a(u,v) &= \int_\Omega \Delta{u} \, \Delta{v}\D\Omega  + \gamma \int_{\Gamma_D} uv \D\Gamma  + \tau \int_{\Gamma_D} \left(\nabla u \cdot \vec n  \right)  \left(\nabla v \cdot \vec n  \right) \D\Gamma +  \int_{\Gamma_D} \Big( u \left( \nabla \left( \Delta v \right) \cdot \vec n \right)   + v \left( \nabla \left( \Delta u \right) \cdot \vec n \right) \Big) \D\Gamma
	\\ 
	& \quad -  \int_{\Gamma_D} \Big( \Delta u \left( \nabla v \cdot \vec n \right) + \Delta v \left(\nabla u \cdot \vec n \right) \Big) \D\Gamma \, ,
	\\ \nonumber
	l(v) & = \int_\Omega v s \D \Omega - \int_{\Gamma_N} v \bar{\lambda} \D\Gamma + \int_{\Gamma_N} \left( \nabla v \cdot \vec n \right) \bar{\kappa}    \D\Gamma + \gamma \int_{\Gamma_D} v \bar{u} \D \Gamma + \tau
	\int_{\Gamma_D} \left( \nabla v \cdot \vec n \right) \bar{t} \D\Gamma 
	\\ 
	& \quad + \int_{\Gamma_D} \left( \nabla \left( \Delta v \right) \cdot \vec n \right) \bar{u} \D\Gamma  - \int_{\Gamma_D} \left(\Delta v \right) \bar{t} \D\Gamma \,.
	\end{align}
\end{subequations}
%
\subsection{Finite element discretisation}
%
We discretise the trial and test functions with the SB-splines as
\begin{equation} \label{eq:femInterpolation}
u^h(\vec x) = \sum_{i=1}^{n_N} N_i( \vec x) \, \alpha_i \quad \text{and} \quad v^h(\vec x) = \sum_{i=1}^{n_N} N_i( \vec x) \, \beta_i \,.
\end{equation}
Introducing the interpolation equation~\eqref{eq:femInterpolation} into the weak form of Poisson equation~\eqref{eq:femPoissonWeakForm} or biharmonic equation~\eqref{eq:femBiharmonicWeakForm} yields a system of linear equations with the unknowns~$\alpha_i$. For instance, the bilinear form~$a(u^h, v^h)$ for the Poisson equation becomes after discretisation
\begin{equation}
a( u^h, v^h ) = \sum_{i=1}^{n_N} \sum_{j=1}^{n_N} 
\alpha_i \left(  \int_\Omega \nabla N_i \nabla N_j \D \Omega \,  + \gamma \int_{\Gamma_D} N_i N_j \D \Gamma -  \int_{\Gamma_D} \Big( N_i \left( \nabla N_j \cdot \vec n \right) + N_j \left(\nabla N_i \cdot \vec n \right) \Big) \D\Gamma  \right) \beta_j \, .
\end{equation}
As usual, the domain integral is evaluated numerically after splitting it into~$n_{el}$ element contributions
\begin{equation}
a( u^h, v^h ) =\sum_{k = 1}^{n_{el}} \left( \sum_{i=1}^{n_N} \sum_{j=1}^{n_N} 
\alpha_i \left(  \int_{\Omega_k} \nabla N_i \nabla N_j \D \Omega_k \,  + \gamma \int_{(\Gamma_D)_k} N_i N_j \D \Gamma_k -  \int_{(\Gamma_D)_k} \Big( N_i \left( \nabla N_j \cdot \vec n \right) + N_j \left(\nabla N_i \cdot \vec n \right) \Big) \D\Gamma_k  \right) \beta_j  \right)\, .
\end{equation}

\section{Mesh refinement \label{sec:refinement}}
\begin{figure}[h]
	\centering
	\subfloat[][Step $1$ of refinement \label{fig:midPoint1}] {
		\includegraphics[width=0.225 \textwidth]{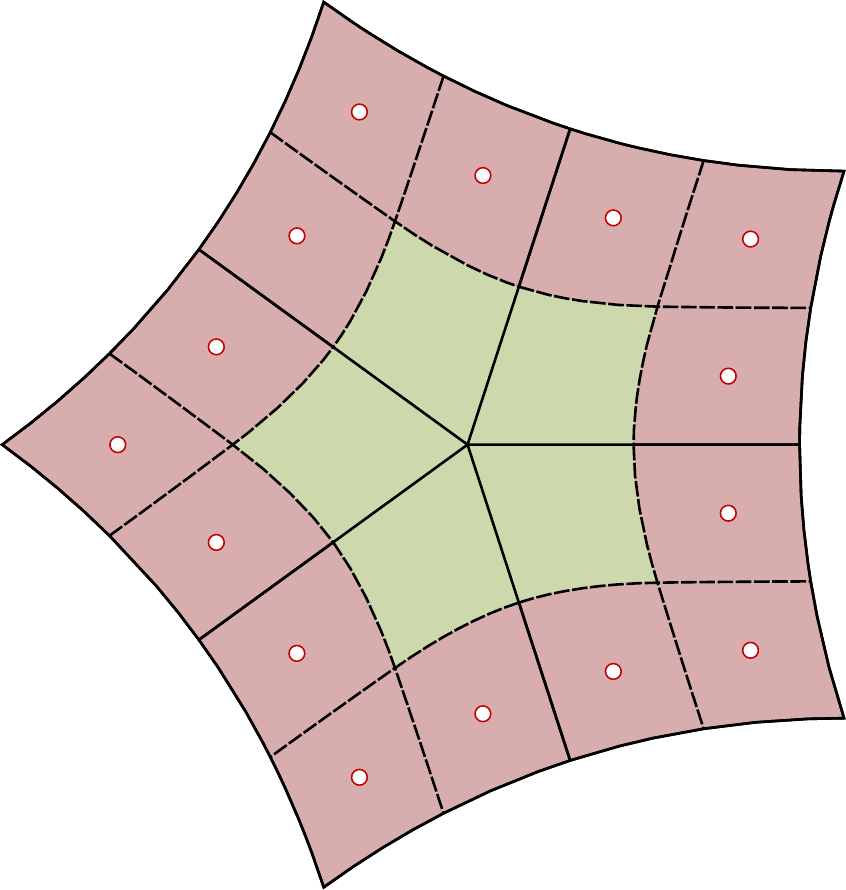} }
	\hfil
	\subfloat[][Step $2$ of refinement \label{fig:midPoint2}] {
		\includegraphics[width=0.225 \textwidth]{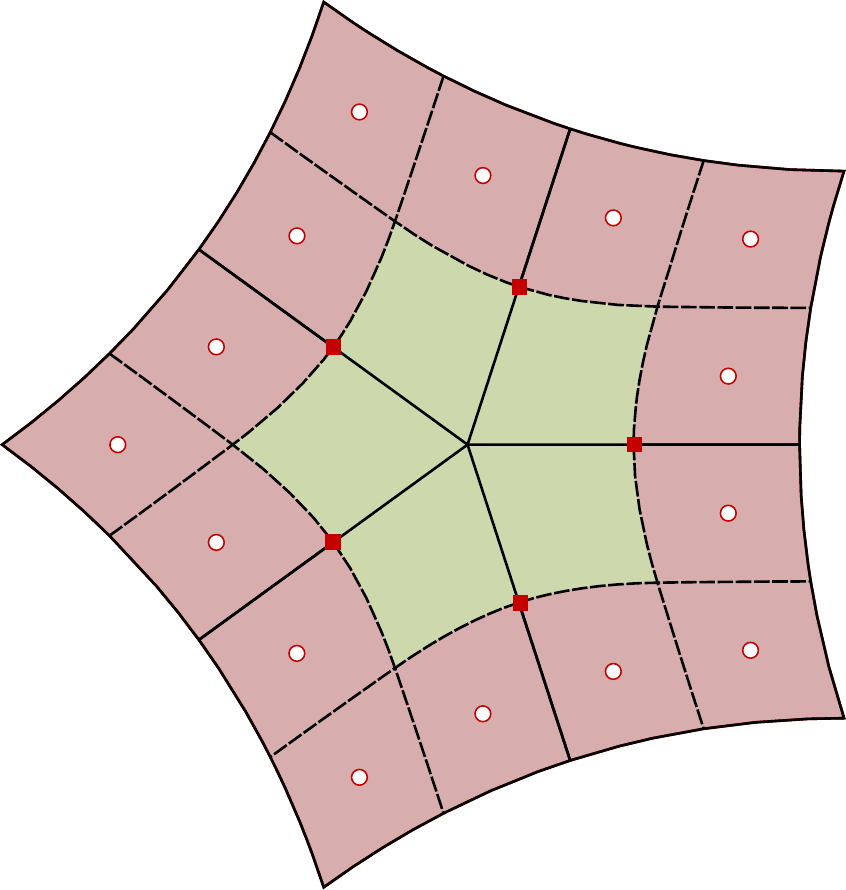} }
	\hfil
	\subfloat[][Step $3$ of refinement \label{fig:midPoint3}] {
		\includegraphics[width=0.225 \textwidth]{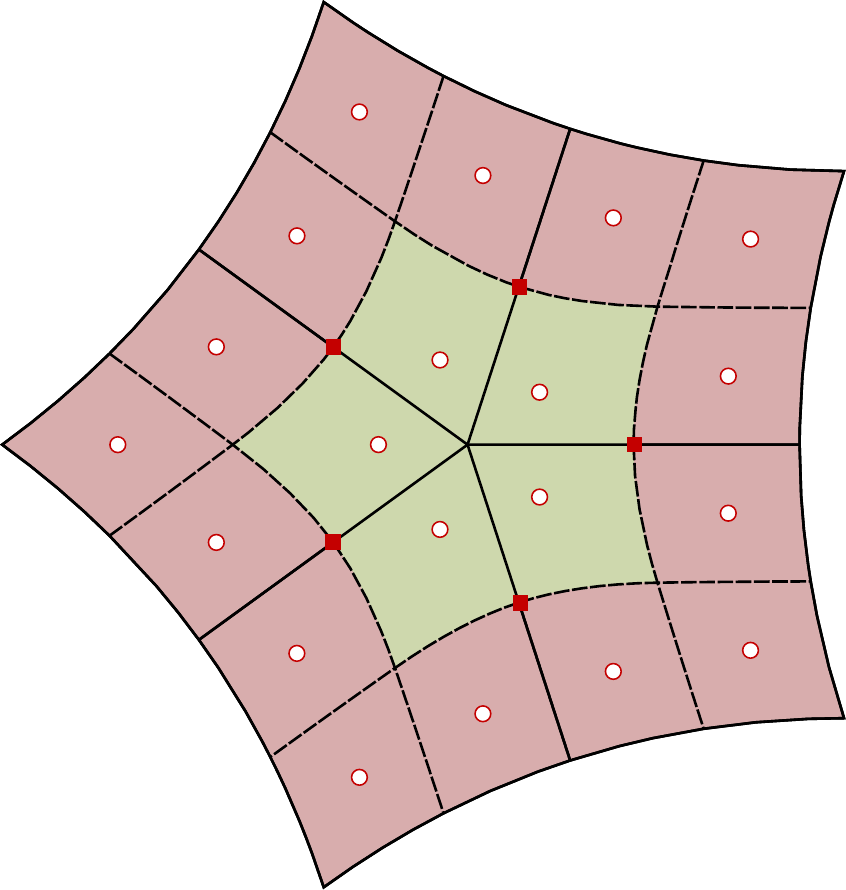} }
	\caption{Refinement of the~$1$-neighbourhood of an extraordinary vertex with valence~$v = 5$. The element boundaries of the coarse and the refined B\'ezier meshes are depicted as solid line and dashed line, respectively. The refinement scheme consists of three steps. First, the $2$-neighbourhood mixed B-spline control vertices of the refined B\'ezier mesh shown in (a) are obtained using tensor-product knot insertion. Second, we compute the midpoints such that the curve lengths of the~$1$-neighbourhood edges of the coarse B\'ezier mesh are approximately bisected shown in (b). Third, the~$1$-neighbourhood mixed B-spline control vertices of the refined B\'ezier mesh shown in (c) are computed such that the midpoints determined in (b) are interpolated.\label{fig:midPoint}}
\end{figure}

%
We use for the unstructured quadrilateral mesh the non-nested refinement scheme by Toshniwal~\cite{toshniwal2022quadratic}. Given a set of mixed B-spline control vertices from the coarse B\'ezier mesh, the objective is to obtain a new set of mixed B-spline control vertices for defining the refined B\'ezier mesh. Away from the~$1$-neighbourhood of an extraordinary vertex, that is, where a tensor product structure is locally present, the knot insertion algorithm is used. The refinement of the $1$-neighbourhood of an extraordinary vertex of valence~$v$ consists of three steps shown in Figure~\ref{fig:midPoint}. First, the~$v \times 3$ mixed B-spline control vertices at the $2$-neighbourhood of the refined B\'ezier mesh are obtained from the knot insertion algorithm as shown in Figure~\ref{fig:midPoint1}. As a result, only the~$v$ mixed B-spline control vertices at the $1$-neighbourhood of the refined B\'ezier mesh remain to be selected. In particular, the remaining~$v$ mixed B-spline control vertices are selected such that the~$v$ midpoints of the edges at the~$1$-neighbourhood of the coarse B\'ezier mesh are interpolated. Therefore, the second step is to estimate the midpoints shown in Figure~\ref{fig:midPoint2} using, for instance, a root-finding algorithm together with a parametrisation for the edge length. Subsequently, in the third step a~$v \times v$ linear system of equations is solved for the~$v$ mixed B-spline control vertices at the $1$-neighbourhood of the refined B\'ezier mesh shown in Figure~\ref{fig:midPoint3}.

The~$v \times v$ linear system is invertible for the case of odd valences~$v = 3, 5,\dotsc$ but not the case of even valences $v = 6, 8,\dotsc$. For even valences, the~$v \times v$ linear system has a rank of~$v - 1$. For the case of even valences, following~\cite{toshniwal2022quadratic} we constrain~$1$ of the~$v$ mixed B-spline control vertices so that the~$v \times v$ linear system has a unique solution. Evidently, the choice of the constrained mixed B-spline control vertex is not arbitrary. For instance, choosing to constraint a mixed B-spline control vertex that is far away from the extraordinary vertex can distort the refined B\'ezier mesh. To avoid any mesh distortion, in this paper, we first select~$1$ of the~$v$ extraordinary elements from the coarse B\'ezier mesh. After that, assuming that the extraordinary vertex is located at the reference element origin~$\vec \eta = (0, 0)$ of the selected extraordinary element, we constrain the mixed B-spline control vertex at $\vec \eta = (\eta_1, \eta_2)$ where~$\eta_1$ and~$\eta_2$ are decided case by case, i.e.\ depending on the coarse B\'ezier mesh. For example, in Section~\ref{sec:exampleVGon}, we observe that the choice of~$\eta_1 = \eta_2$ with~$0.125 \leq \eta_1 \leq 0.25$ generally preserves the mesh quality after refinement.
%
\section{Arbitrary joint and prism valences \label{sec:joint}}
\begin{figure}[h!]
	\centering
	\subfloat[][Hexahedral mesh \label{fig:jointA}] {
		\includegraphics[width=0.285 \textwidth]{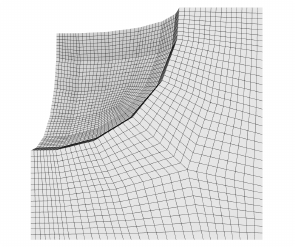} } 
	\hfil
	\subfloat[][Extraordinary edges and vertices \label{fig:jointB}] {
		\includegraphics[width=0.285\textwidth]{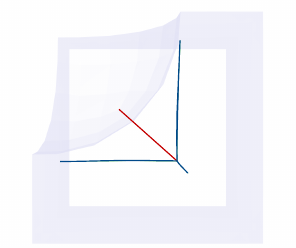} }
	\hfil
	\subfloat[][Extraordinary hexahedra \label{fig:jointC}] {
		\includegraphics[width=0.285\textwidth]{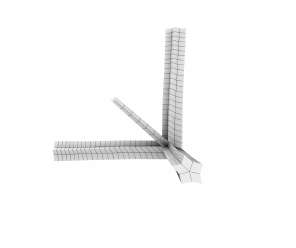} }
	\hfil
	\subfloat[][Extraordinary prisms \label{fig:jointD}] {
		\includegraphics[width=0.285\textwidth]{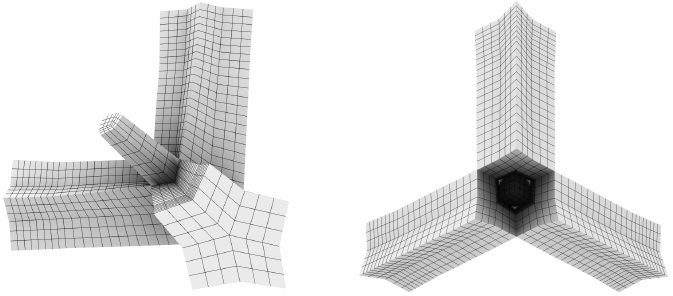} }
	\hfil
	\subfloat[][Extraordinary joint and the attached prisms \label{fig:jointE}] {
		\includegraphics[width=0.4\textwidth]{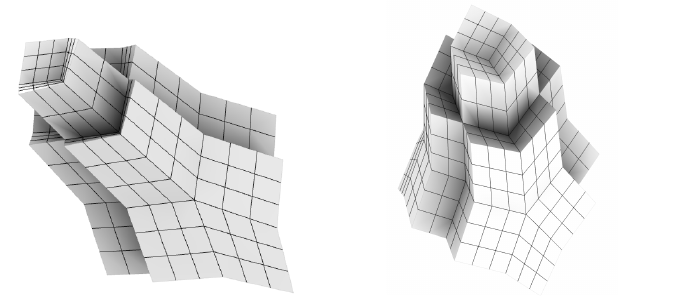} }
	\caption{Truncated box domain with a spherical cavity and (a) its discretisation with an unstructured hexahedral mesh, (b) extraordinary edges and vertices of the mesh, (c) extraordinary hexahedra, (d) extraordinary prisms and (e) extraordinary joint and the attached prisms. In (b) the four sets of extraordinary edges of valence~$e = 3$ and~$e = 5$ are coloured in red and blue, respectively.   \label{fig:joint}}
\end{figure}

%
We briefly demonstrate that the approach discussed in Section~\ref{sec:weightJoint} applies to an extraordinary joint with arbitrary edge valences~$e$ and vertex valence~$v$. As a concrete example, Figure~\ref{fig:joint} illustrates the construction of weight functions for an unstructured hexahedral mesh of a truncated box domain with a spherical cavity. The unstructured hexahedral mesh shown in~\ref{fig:jointA} consists of one set of extraordinary edges of valence~$e = 3$ and three sets of extraordinary edges of valence~$e = 5$ shown in~\ref{fig:jointB}. The valences of the extraordinary edges can be verified from the extraordinary hexahedra of the mesh shown in~\ref{fig:jointC}. Overall, the unstructured hexahedral mesh consists of an extraordinary joint where the extraordinary edges of valence~$e = 3$ and~$e = 5$ meet at a vertex of valence~$v = 10$. As discussed in Section~\ref{sec:weightJoint}, we first require that the support of the prism weight functions do not overlap at the extraordinary joint, see Figure~\ref{fig:jointD} depicting the union of the prism weight function supports. Subsequently, the extraordinary joint weight function is defined over the set of hexahedra in Figure~\ref{fig:jointE}. Conceptually, the extraordinary joint weight function defined over the set of hexahedra in Figure~\ref{fig:jointE} resembles that shown in Figure~\ref{fig:threeDTetWeightwQ1}.

\bibliographystyle{elsarticle-num-names}
\bibliography{psplines}

\end{document}